\documentclass[final]{elsarticle}

\usepackage[dvipsnames]{xcolor}

\usepackage[utf8]{inputenc}
\usepackage[english]{babel}
\usepackage{csquotes}
\usepackage{lmodern}
\usepackage[T1]{fontenc}
\usepackage[final]{microtype}
\usepackage[mathscr]{euscript} 

\usepackage{amsmath}
\usepackage{amssymb}
\usepackage{amsthm}
\usepackage{mathtools}
\usepackage{stmaryrd}
\usepackage{subcaption}
\SetSymbolFont{stmry}{bold}{U}{stmry}{m}{n}

\usepackage[colorinlistoftodos,prependcaption,textsize=small]{todonotes}

\usepackage{pgf}
\usepackage{tikz}
\usepackage{tikz-cd}
\usetikzlibrary{arrows, matrix}
\usepackage{subcaption}
\usepackage{multicol}
\usepackage{multirow}
\usepackage{pgfplots}

\usepackage{caption}
\usepackage[position=b]{subcaption}
\usepackage[linesnumbered,ruled,vlined]{algorithm2e}

\theoremstyle{definition}
\newtheorem{definition}{Definition}[section]
\theoremstyle{plain}

\newtheorem{lemma}[definition]{Lemma}

\newtheorem{corollary}[definition]{Corollary}
\theoremstyle{remark}
\newtheorem{remark}[definition]{Remark}

\newcommand{\bb}{\pmb}

\newcommand{\R}{{\mathbb R}}
\newcommand{\C}{{\mathbb C}}
\newcommand{\K}{{\mathbb K}}
\newcommand{\norm}[1]{{\left\lVert #1 \right\rVert}}
\newcommand{\abs}[1]{{\left\lvert #1 \right\rvert}}

\newcommand{\opd}{\,\operatorname{d}}

\renewcommand{\S}{{\mathbb S}}

\renewcommand{\hat}{\widehat}

\usetikzlibrary{patterns}
\usepackage{microtype}
\usepackage{placeins}

\title{A Fast Isogeometric BEM for the Three Dimensional Laplace- and Helmholtz Problems}

\author{
	  J\"urgen D\"olz$^{1}$,  Helmut Harbrecht$^{1}$,  Stefan Kurz$^{2}$, \\ Sebastian Sch\"ops$^{2}$,  Felix Wolf$^{~2}$\\
\small$^{1}$ Universit\"at Basel,
Departement Mathematik und Informatik\\
\small$^{2}$ Technische Universit\"at Darmstadt,
Institute TEMF \& Graduate School CE
}
\date{\today}

\begin{document}

\begin{abstract}
	We present an indirect higher order boundary element method utilising NURBS mappings for exact geometry representation and an interpolation-based fast multipole method for compression and reduction of computational complexity, to counteract the problems arising due to the dense matrices produced by boundary element methods.
	By solving Laplace and Helmholtz problems via a single layer approach we show, through a series of numerical examples suitable for easy comparison with other numerical schemes, that one can indeed achieve extremely high rates of convergence of the pointwise potential through the utilisation of higher order B-spline-based ansatz functions.
\end{abstract}

\maketitle


\section{Introduction}

In the search of a method incorporating simulation techniques into the design workflow of industrial development, \cite{Hughes_2005aa} proposed the concept of \emph{Isogeometric Analysis (IGA)} to unite \emph{Computer Aided Design (CAD)} and \emph{Finite Element Analysis (FEA)}; making it possible to perform simulations directly on geometries described by volumetric NURBS parametrisations.

Nontheless, many CAD systems use boundary representations only, thus volumetric parametrisations often have to be constructed solely for the purpose of simulation \cite{CornoSomehowitisimpossibleformetofindthisreferenceinthedatabasebibfile}.
The boundary parametrisation, however, can be easily exported from CAD. Thus, an approach via isogeometric boundary element methods seems natural. 
Indeed, isogeometric boundary element methods have recently been tested successfully for multiple engineering applications, see e.g.~\cite{BEER2017418,MZB15,Simpson_2012aa,2017arXiv170407128S}.

\begin{figure}[h]
	\centering
	\begin{subfigure}{\textwidth}\centering
		\includegraphics[width=.7\textwidth]{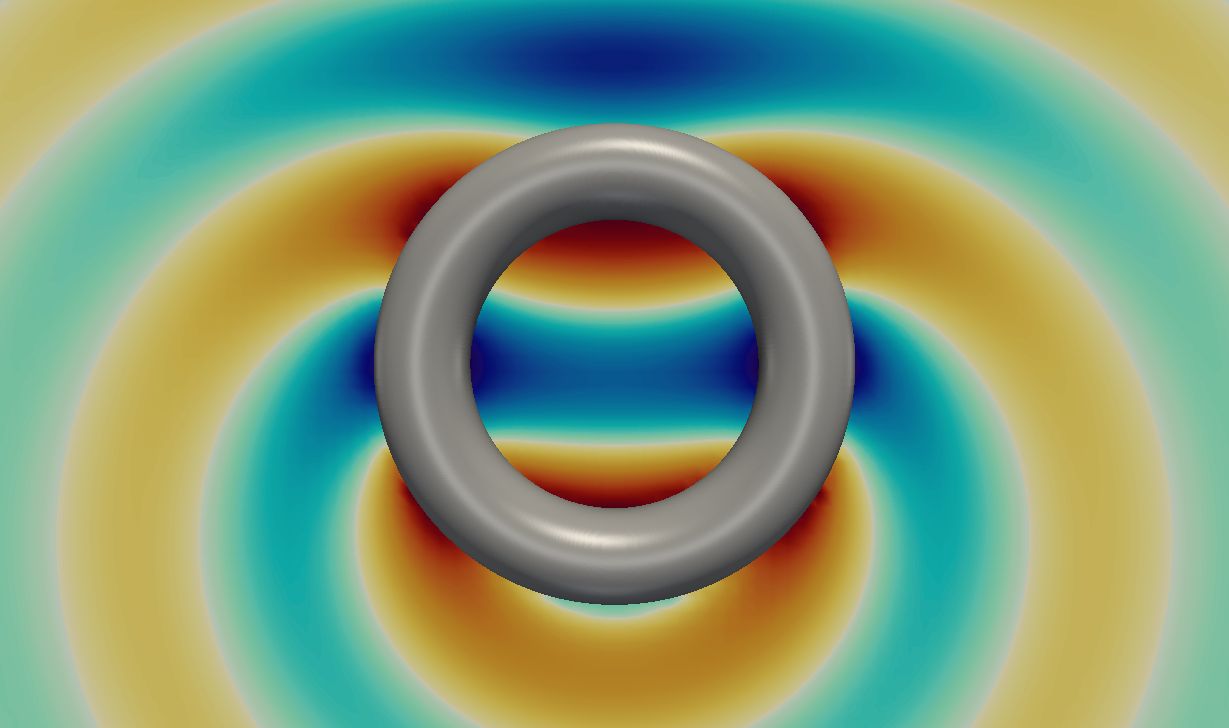}
		\caption{Scattered wave.}
	\end{subfigure}
	\begin{subfigure}{\textwidth}\centering
		\includegraphics[width=.7\textwidth]{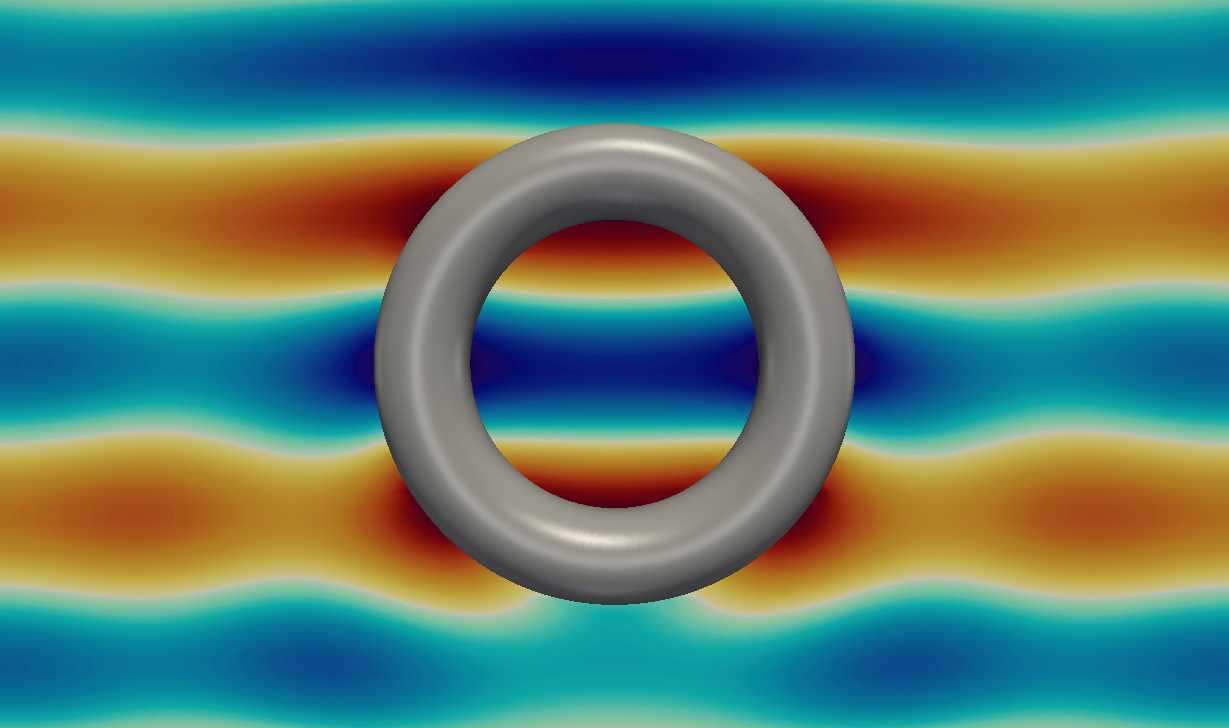}
		\caption{Scattered wave with incident wave.}
	\end{subfigure}
	\caption{Real part of the exterior potential of a plane wave Helmholtz scattering problem computed with the method discussed in this paper.}\label{fig::scattering}
\end{figure}

The utilisation of parametric mappings in numerical implementations of boundary element methods is not new, going back further than the introduction of the isogeometric concept \cite{harbrech2001}.
Parametric mappings avoid the problem of a slow convergence of the geometry due to the limited polynomial approximation of the geometry \cite{Weggler_2011aa}. Thus, they encourage the application of higher order Galerkin schemes.

Through the parametric mappings, a tensor product structure on the geometries is induced, making it possible to define patchwise tensor product B-spline bases of high order and regularity \cite{Schumaker_1981aa}.

One of the major downsides of the application of boundary element methods is, that the integral operators involved yield dense discrete systems.
To counteract the dense matrices, so-called \emph{fast methods} must be utilised for compression and efficiency. As shown in \cite{doelz2016,Harbrecht_2013aa}, the tensor product structure induced by the mappings can be exploited to achieve an efficient implementation of compression techniques such as $\mathcal H$-Matrices or the Fast Multipole Method \cite{greengard2009fast,greengard1987fast,Kurz_2007aa}. Thus, an isogeometric boundary element method promises runtimes which can compete with classical discretisation methods.

\begin{figure}
	\centering
	\includegraphics[width=.7\textwidth]{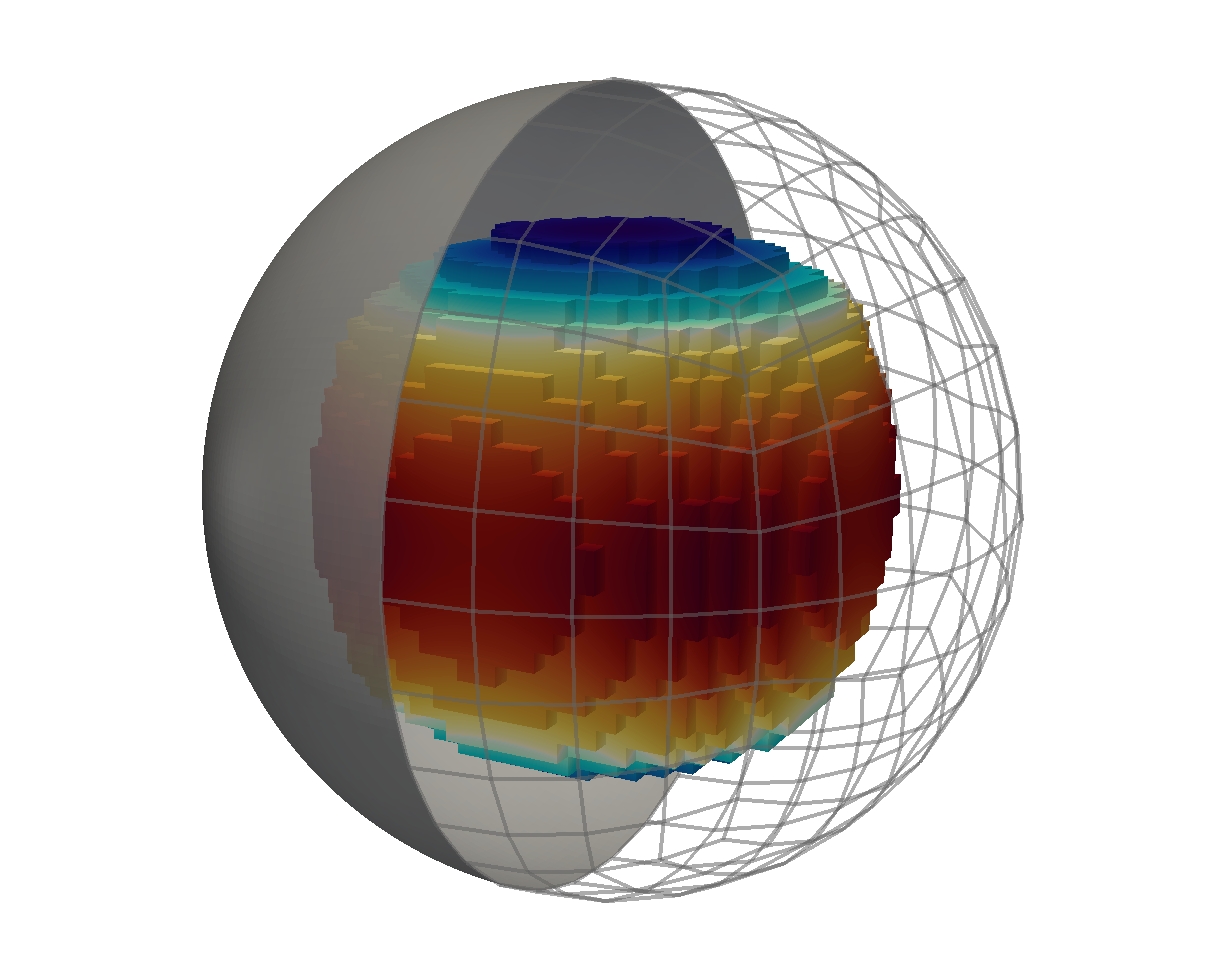}
	\caption{Potential evaluation of one of the Laplace test cases.}
	\label{fig::spherepic}
\end{figure}

However, both analysis and application of higher order B-splines (in view of regularity requirements, convergence of the potential within the domain) and multipole compression techniques are not established. 

Thus, in this paper, we want to discuss numerical examples computed via an indirect parametric boundary element method utilising ansatz functions of higher polynomial degrees for both the Laplace and Helmholtz problem in three dimensions on different geometries, cf.~Figure \ref{fig::scattering}, paired with an interpolation based fast multipole method developed in \cite{doelz2016}.

In Section \ref{sec::IGA}, we will review the basic concepts of isogeometric analysis, defining the B-spline basis and NURBS mappings for the geometry. Utilising the transformations induced by the geometry mappings, we will then define the discrete function space used in the numerical experiments.

Section \ref{sec::analysis} then introduces the analytic setting, briefly introducing the approach to Helmholtz and Laplace problems via indirect boundary element methods, and explaining how to recast the arising variational problem w.r.t.~the reference domain of the isogeometric mappings of the previous section.

Afterwards, in Section \ref{sec::impl}, we will discuss the major novelties of the method utilised in this paper and in \cite{doelz2016}, which include an interpolation based multipole method explained in Section \ref{subsec::fmm} and an approach via a discontinuous superspace as presented in Section \ref{subsec::embedding}.

Following the explanation of the method utilised for the numerical solution of the problems, we will explain and discuss the numerical test cases in Section \ref{sec::numex}, where the data for the Laplace test cases can be found in Section \ref{subsec::laplnumex}. Section \ref{subsec::perturbed} repeats some examples of the Laplace case with perturbed mappings. The Helmholtz results are then presented in Section \ref{subsec::helmnumex}.

Finally, we will conclude our findings in Section \ref{sec::final}.

\FloatBarrier

\section{Isogeometric Analysis}\label{sec::IGA}

What follows is a brief introduction to the basic concepts of isogeometric analysis, starting with the definition of the B-spline basis, followed by the description of the geometry using NURBS.

\subsection{B-splines}
For this section, let $\K$ be either $\R$ or $\C$. The original definitions (or equivalent notions) and proofs, as well as basic algorithms, can be found in most of the standard spline- and isogeometric literature \cite{Cottrell_2009aa,Hughes_2005aa,Lee_1996aa,Piegl_1997aa,Schumaker_1981aa}.

\begin{figure}
	\centering
	\begin{subfigure}{.49\textwidth}
		\begin{tikzpicture}
			\begin{axis}[
					xmin = -.1,
					xmax = 1.1,
					ymin = -.1,
					ymax = 1.1,
					width=1\columnwidth,
					height=.66\columnwidth,
					grid=major,
					legend style={
							at={(.5,1)},
							anchor=south}
				]
				\addplot[red,ultra thick,mark=none,domain = 0:1/3] {1};
				\addplot[teal,ultra thick,mark=none,domain = 1/3:2/3] {1};
				\addplot[blue,ultra thick,mark=none,domain = 2/3:1] {1};
			\end{axis}
		\end{tikzpicture}
		\caption{$p=0$, $\Xi=[0,1/3,2/3,1]$.}
	\end{subfigure}
	\begin{subfigure}{.49\textwidth}
		\begin{tikzpicture}
			\begin{axis}[
					xmin = -.1,
					xmax = 1.1,
					ymin = -.1,
					ymax = 1.1,
					width=1\columnwidth,
					height=.66\columnwidth,
					grid=major,
					legend style={
							at={(.5,1)},
							anchor=south}
				]
				\addplot[red,ultra thick,mark=none,domain = 0:1/3] {3*-(x-1/3)};
				\addplot[teal,ultra thick,mark=none,domain = 0:1/3] {3*(x)};
				\addplot[blue,ultra thick,mark=none,domain = 1/3:2/3] {3*(x-1/3)};
				\addplot[orange,ultra thick,mark=none,domain = 2/3:1]{3*(x-2/3)};
				\addplot[teal,ultra thick,mark=none,domain = 1/3:2/3] {(1-3*x)+1};
				\addplot[blue,ultra thick,mark=none,domain = 2/3:1] {3*(1-x-1/3)+1};
			\end{axis}
		\end{tikzpicture}
		\caption{$p=1$, $\Xi=[0,0,1/3,2/3,1,1]$.}
	\end{subfigure}\\[.5cm]\centering
	\begin{subfigure}{\textwidth}
		\begin{tikzpicture}[scale = .95]
			\begin{axis}[
					xmin = -.1,
					xmax = 1.1,
					ymin = -.1,
					ymax = 1.1,
					width=1\columnwidth,
					height=.3397\columnwidth,
					grid=major,
					legend style={
							at={(1,.5)},
							anchor=west}
				]
				\addplot[red,ultra thick,mark=none,domain = 0:1/3]{(3*(x)-1)^2} ;
				\addplot[teal,ultra thick,mark=none,domain = 0:1/3]{2*(3*(x))*(1-3*x)+.5*(3*x)^2};
				\addplot[blue,ultra thick,mark=none,domain = 0:1/3]{.5*(x*3)^2};
				\addplot[orange,ultra thick,mark=none,domain = 2/3:1]{2*(3*((x-2/3)))*(1-3*(x-2/3))+.5*(1-3*(x-2/3))^2};
				\addplot[brown,ultra thick,mark=none,domain = 2/3:1]{(3*(x-2/3))^2};
				\addplot[teal,ultra thick,mark=none,domain = 1/3:2/3]{.5*(1-3*(x-1/3))^2 };
				\addplot[blue,ultra thick,mark=none,domain = 1/3:2/3]{-(((x*3)-1)-1)*((x*3)-1)+.5};
				\addplot[blue,ultra thick,mark=none,domain = 2/3:1]{.5*(3-(x*3))^2};
				\addplot[orange,ultra thick,mark=none,domain = 1/3:2/3]{.5*(3*(x-1/3))^2 };
			\end{axis}
		\end{tikzpicture}
		\caption{$p=2$, $\Xi=[0,0,0,1/3,2/3,1,1,1]$.}
	\end{subfigure}
	\caption{B-spline bases for $p=0,1,2$ and open knot vectors with interior knots $1/3$ and $2/3$.}\label{fig::splines}
\end{figure}
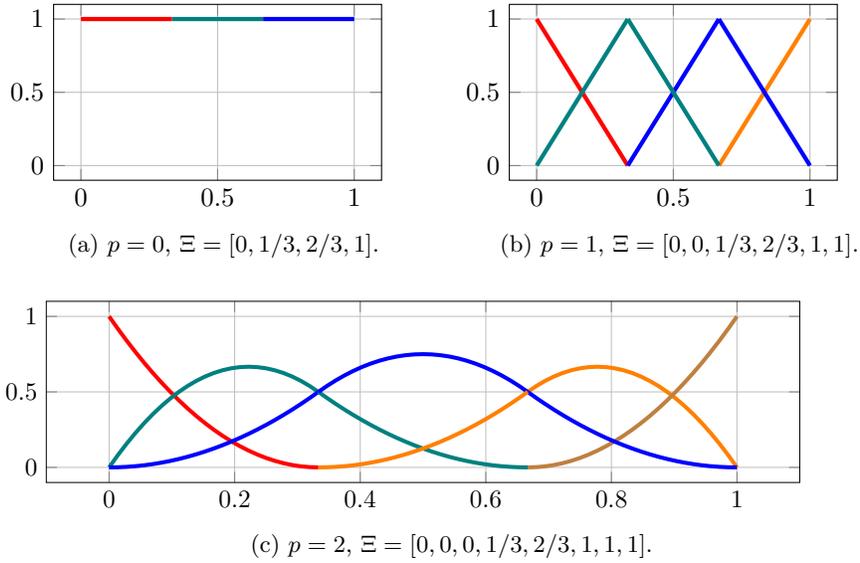

\begin{definition}
	Let $0\leq p\leq k$. We define a \emph{$p$-open knot vector} as a set
	\begin{align}
		\Xi = \big[\underbrace{\xi_0 = \cdots =\xi_{p}}_{=0}\leq \cdots \leq \underbrace{\xi_{k}=\cdots =\xi_{k+p}}_{=1}\big]\in[0,1]^{k+p+1},
	\end{align}
	where $k$ denotes the number of control points.

	We can then define the basis functions $\left(b_j^p\right)_{0\leq j< k}$ for $p=0$ as
	\begin{align}
		b_j^0(x) & =\begin{cases}
			1, & \text{if }\xi_j\leq x<\xi_{j+1}, \\
			0, & \text{otherwise,}
		\end{cases}
		\intertext{ and for $p>0$ via the recursive relationship}
		b_j^p(x) & = \frac{x-\xi_j}{\xi_{j+p}-\xi_j}b_j^{p-1}(x) +\frac{\xi_{j+p+1}-x}{\xi_{j+p+1}-\xi_{j+1}}b_{j+1}^{p-1}(x),
	\end{align}
	cf.~Figure \ref{fig::splines}. A \emph{B-spline} is then defined as a function
	$$f(x) = \sum_{0\leq j< k}p_jb_j^p(x),$$
	where $\lbrace p_j\rbrace_{0\leq j< k}\subset\K$ denotes the set of \emph{control points.} If one sets  $\lbrace p_j\rbrace_{0\leq j< k}\subset\K^d$, then $f$ will be called a \emph{B-spline curve}.
\end{definition}
\begin{definition}
	Let $\Xi$ be a $p$-open knot vector containing $k+p+1$ elements. We define the \emph{B-spline space} $S_{p}(\Xi)$ as the space spanned by $(b_j^p)_{0\leq j< k}$.
\end{definition}

\begin{definition}
	For a knot vector $\Xi,$ let $h_j\coloneqq  \xi_{j+1}-\xi_{j}.$ We define the \emph{mesh size} $h$ to be the maximal distance $h\coloneqq \max_{p\leq j < k}h_j$ between neighbouring knots.
	We call a knot vector \emph{quasi uniform}, when there exists a constant $\theta\geq 1$ such that for all $p\leq j < k$ the ratio $h_j\cdot h_{j+1}^{-1}$ satisfies $\theta^{-1}\leq h_j\cdot h_{j+1}^{-1} \leq \theta.$
\end{definition}

B-splines on higher dimensional domains are constructed through simple tensor product relationships for $\bb p_{j_1,\dots j_\ell} \in \K^n$ via
\begin{align}
	\bb f(x_1,\dots ,x_\ell)=\sum_{0\leq j_1< k_1}\dots\sum_{0\leq j_\ell< k_\ell} \bb p_{j_1,\dots,j_\ell} \cdot b_{j_1}^{p_1}(x_1)\cdots b_{j_\ell}^{p_\ell}(x_\ell),\label{def::tpspline}
\end{align}
which allows \emph{tensor product B-spline spaces} to be defined as $$S_{p_1,\dots,p_\ell}(\Xi_1,\dots,\Xi_\ell).$$

Throughout this paper, we will reserve the letter $h$ for the mesh size. All knot vectors will be assumed to be $p$-open and quasi uniform, such that the usual spline theory is applicable \cite{BuffaActa,Piegl_1997aa,Schumaker_1981aa}.

\subsection{Description of the Geometry and Discrete Space}

Regarding the description of the geometry and an adequate choice of discrete function spaces, we review the most important properties required within this paper.

The geometry will be assumed to be closed and Lipschitz continuous. For the remainder of this article, we assume that the boundary $\Gamma=\bigcup_{j<n}\Gamma_j$ is given patchwise, i.e.~that it is induced by (smooth, nonsingular, bijective) diffeomorphisms
\begin{align}
	\bb F_j\colon \widehat \Omega = [0,1]^2 \to \Gamma_j \subset \R^3.\label{def::geom}
\end{align}
In the spirit of isogemetric analysis, these mappings will be given by NURBS mappings, i.e.,~by
\begin{align}
	\bb F_j(x,y)\coloneqq \sum_{0\leq j_1<k_1}\sum_{0\leq j_2<k_2}\frac{\bb c_{j_1,j_2} b_{j_1}^{p_1}(x) b_{j_2}^{p_2}(y) w_{j_1,j_2}}{ \sum_{i_1=0}^{k_1-1}\sum_{i_2=0}^{k_2-1} b_{i_1}^{p_1}(x) b_{i_2}^{p_2}(y) w_{i_1,i_2}},
\end{align}
for control points $\bb c_{j_1,j_2}\in \R^3$ and weights $w_{i_1,i_2}>0.$
We will moreover require that, for any interface $D = \Gamma_j\cap \Gamma_i \neq \emptyset$, the NURBS mappings coincide, i.e.~that, up to rotation of the reference domain, one finds $\bb F_j(\cdot,1) \equiv \bb F_i(\cdot,0).$

The mappings of \eqref{def::geom} give rise to the transformations
\begin{align}
	\iota_j(f){}\coloneqq{} & {}f\circ  \bb F_j{},
\end{align}
which can be utilised to define discrete spaces patchwise, by mapping the space of tensor product B-splines as in \eqref{def::tpspline} with
\begin{align}
	\Xi_{p,m} \coloneqq \big[ \underbrace{0,\dots,0}_{p+1\text{ times}},{1}/{2^m},\dots,  {(2^m-1)}/{2^m},\underbrace{1,\dots,1}_{p+1\text{ times}}\big]
\end{align}
to the geometry. Here the variable $m$ denotes the level of uniform refinement. For our purposes, the global function space on $\Gamma$ can thus be defined by
\begin{align}
	\mathbb S_{p,m}(\Gamma) \coloneqq \left\lbrace f\in H^{-1/2}(\Gamma)\colon f|_{\Gamma_j} \equiv \iota_j^{-1}(g)\text{ for some }g\in S_{p,p}(\Xi_{p,m},\Xi_{p,m})\right\rbrace,\label{def::space}
\end{align}
as commonly done in the isogeometric literature, see e.g. \cite{Buffa_2011aa,Buffa_2010aa,Buffa_2013aa}. Note that the spline space $\mathbb S_{p,m}(\Gamma)$ is of dimension $n\cdot(2^m + p )^2$, where $n$ denotes the number of patches involved in the description of the geometry.

\section{Analytic Setting} \label{sec::analysis}

The relevant Sobolev spaces will be introduced briefly. Let $\Omega\subseteq \R^d$ be some compact Lipschitz domain. Let $Df$ denote the weak derivative of $f$. Following conventions as in \cite{Buffa_2003ab} or \cite{Nezza_2012aa}, we will set $H(\Omega)=H^0(\Omega)= L^2(\Omega)$, and, for any integer $m,$ we will set $H^m(\Omega)=\lbrace f\in L^2(\Omega)\colon D f \in H^{m-1}(\Omega)\rbrace$. Since these definitions have to be considered standard, we will only refer to the established literature \cite{Adams_1978aa,Werner_2011aa}. 

For non integer valued $s>0,$ one can define the corresponding Sobolev spaces via interpolation \cite{Bergh_1976aa} or, equivalently, via the Sobolev-Slobodeckij semi-norm \cite{McLean_2000aa}.

On geometries $\Gamma$ given by families of mappings as considered by \eqref{def::geom}, let $H^{-s}(\Gamma)$ denote the dual space of $H^s(\Gamma)$ for any $s\geq 0$.

\subsection{Indirect Boundary Element Methods}\label{subsec::indirect}

We are interested in the function spaces on some Lipschitz continuous boundary $\Gamma = \partial \Omega$. Throughout this paper, we will assume that there exists a patchwise parametrisation of $\Gamma$ as in \eqref{def::geom}.

Let $\gamma_0$ denote the trace operator, which is known to be linear, continuous and surjective with a continuous right inverse for $\gamma_0\colon H^s(\Omega)\to H^{s-1/2}(\Gamma)$ and $1/2<s<3/2$ \cite{McLean_2000aa}.

Consider one of the following Dirichlet problems. Let $\Omega$ be a compact domain with boundary $\Gamma$. Given a function $g$ on $\Gamma,$ find a function $u$ on $\Omega$ such that, in case of the \emph{interior Laplace problem},
\begin{align}
	\begin{aligned}
		\Delta u  & = 0\qquad\text{in }\Omega,\\
		\gamma_0u & = g\qquad\text{on }\Gamma,\label{problem::laplace}
	\end{aligned}
\end{align}
and, in case of the \emph{exterior Helmholtz problem},
\begin{align}
	\begin{aligned}
		\Delta u + \kappa^2 u & = 0\qquad\text{in }\mathbb R^3\setminus \Omega,\\
		\gamma_0u            & = g\qquad\text{on }\Gamma,\label{problem::helmholtz}
	\end{aligned}
\end{align}
together with the Sommerfeld radiation condition, cf.~\cite{sauter_2010aa}, are fulfilled, where $\kappa$ denotes the wavenumber.
One can derive the well known indirect representation formulae in the form of
\begin{align}
	u(\bb x) = \tilde V(w)(\bb x) \qquad \text{on }\R^3\setminus \Gamma,\label{eq::repform}
\end{align}
for some unknown density $w\in H^{-1/2}(\Gamma)$, where the single layer potential operator $\tilde V$ is defined as
\begin{align}
	\tilde V(\mu)(\bb x) & \coloneqq	\int_\Gamma  u^*(\bb x,\bb y) \mu(\bb y) \opd \bb y,
\end{align}
for all $\mu\in H^{-1/2}(\Gamma)$, and $u^*(\bb x,\bb y)$ represents fundamental solution
\begin{align}
	u^*(\bb x,\bb y) & \coloneqq \frac{1}{4\pi\norm{\bb x-\bb y}}  \label{eq::fundament1}
	\intertext{for the Laplace problem \eqref{problem::laplace}, and }
	u^*(\bb x,\bb y) & \coloneqq \frac{\exp(i\kappa\norm{\bb x-\bb y})}{4\pi\norm{\bb x-\bb y}} \label{eq::fundament2}
\end{align}
for the Helmholtz problem \eqref{problem::helmholtz}, respectively \cite{Rjasanow_2007aa,sauter_2010aa,steinbach_2013aa}. Setting $V=\gamma_0\circ\tilde V,$ the unknown density $w$ can be obtained by solving the variational problem of finding a $\mu\in H^{-1/2}(\Gamma)$ such that
\begin{align}
	\left\langle V(\mu),\nu \right\rangle_\Gamma = \left\langle g,\nu\right\rangle_\Gamma,\qquad \forall \nu \in H^{-1/2}(\Gamma), \label{eq::varform}
	\intertext{leading to the discrete formulation of finding a $\mu_h\in\S_{p,m}(\Gamma)$ such that}
	\left\langle V(\mu_h),\nu_h \right\rangle_\Gamma = \left\langle g,\nu_h\right\rangle_\Gamma,\qquad \forall \nu \in \S_{p,m}(\Gamma). \label{eq::varformdisc}
\end{align}
Existence and uniqueness of the solution of these problems for both, the Laplace and the Helmholtz case, are discussed extensively in the cited literature (for the Helmholtz case under the assumption that $\kappa$ is not a resonant wavenumber, cf.~\cite{Engleder_2006aa}).

The indirect approach has two main advantages: It requires one boundary operator only. Thus, a second operator needs not to be evaluated in an implementation. The other advantage is the ability to solve the interior and exterior problem simultaneously, since the unknown quantity corresponds to the jump relation of the interior and exterior density. However, with this comes a drawback: Even if one of the respective problems, either interior or exterior, is regular, the quality of the solution will suffer if the solution in the other region is irregular \cite{steinbach_2013aa}. Thus, one can expect a quality impact on problems with, e.g., a nonsmooth domain. \label{section::bemapproach}

\begin{remark}
	The insensitivity of the indirect approach to the choice of interior or exterior domain plays well with the idea of isogeometric analysis. 
	Within a design workflow, no process to distinguish exterior and interior domain needs to take place, merely evaluation points need to be chosen. 
	Moreover, since the right hand side does not involve an integral operator as in a direct approach \cite{sauter_2010aa}, for a change of data no integral operators need to be evaluated.
\end{remark}

\subsection{Formulation of the Problem on the Reference Domain}

Defining the \emph{surface measure} of a mapping $\bb F_j$ for $\hat{\bb x} = (x,y)\in [0,1]^2$ as
\begin{align}
	k_j (\hat{\bb x})\coloneqq \norm{\partial_{x}\bb F_j(\hat{\bb x})\times \partial_{y}\bb F_j(\hat{\bb x})},
\end{align}
one can recast the discrete variational formulation \eqref{eq::varformdisc} in terms of the reference domain, since
\begin{align}
  & \left\langle V(\mu_h),\nu_h \right\rangle_\Gamma\nonumber\\
  &\quad= \sum_{j< n} \left\langle V(\mu_h),\nu_h \right\rangle_{\Gamma_{j}}\nonumber\\
  &\quad= \sum_{i,j< n} \int_{\Gamma_i}\int_{\Gamma_j}  u^*(\bb x,\bb y) \mu_h(\bb x)\nu_h(\bb y) \opd \bb y \opd \bb x\\
  &\quad= \sum_{i,j<n} \int_{[0,1]^2}\int_{[0,1]^2}  k_{j}(\hat{\bb x})k_{i}(\hat{\bb y})u^*\Big(\bb F_j(\hat{\bb x}),\bb F_i(\hat{\bb y})\Big) \mu_h\Big(\bb F_j(\hat{\bb x})\Big) \nu_h\Big(\bb F_i(\hat{\bb y})\Big) \opd \hat{\bb y} \opd \hat{\bb x}.\nonumber
\end{align}
Applying a similar reasoning to the right hand side yields
\begin{align}
	\langle g,\nu_h\rangle_\Gamma = \sum_{i<n} \int_{[0,1]^2} k_{i}(\hat{\bb x})g(\bb F_i(\hat{\bb x}))\nu_h(\bb F_i(\hat{\bb x})) \opd \hat{\bb x},
\end{align}
and, thus, the method can operate on the reference domain only.
\subsection{Convergence Estimates of Optimal Order}
\label{subsec::convergence}

Since we aim for a discretisation of $H^{-1/2}(\Gamma),$ matching conditions across patch boundaries are not required, and thus we can prove the following for multipatch domains, based on the single patch results given in \cite{BuffaActa}.

\begin{lemma}\label{cor::Approxcor}

	For a boundary $\Gamma$ given as by \eqref{def::geom}, we define the norms
	\begin{align}
		\norm{f}_{\tilde H^s(\Gamma)} \coloneqq \sum_{j< n} \norm{f|_{\Gamma_j}}_{H^s(\Gamma_j)}
	\end{align}
	and the spaces $\tilde H^s(\Gamma)$ as
	\begin{align}
		\tilde H^s(\Gamma)\coloneqq \left\lbrace f\in L^2(\Gamma) \colon \norm{f}_{\tilde H^s(\Gamma)}<\infty\right\rbrace.
	\end{align}
	Set $s \coloneqq \min_{j\leq n}\max\lbrace{s>0\colon u\in \tilde H^s(\Gamma)\rbrace}.$
	Then, there exists a projection $$\Pi \colon \tilde H^s(\Gamma)\to \S_{p,m}(\Gamma)$$ such that for $p+1\leq s$ one finds
	\begin{align}
		\norm{u -\Pi u}_{H^r(\Gamma)}
		{} & {}\leq  Ch^{p+1-r}\norm{u}_{\tilde H^{s}(\Gamma)},\qquad r \in \lbrace{-1/2,0\rbrace},\label{eq::estimate}
	\end{align}
	where $C$ is a constant depending on $s$, the geometry mapping and the uniformity coefficient of the knot vector.
	\begin{proof}
		By the assumption that $p+1\leq s$ and \cite[Corollary 5.12]{BuffaActa}, we know there exist projectors $\Pi \colon H^s(\Gamma_j)\to \S_{p,m}(\Gamma_j)$ such that
		\begin{align}
			\norm{u -\Pi  u}_{H^r(\Gamma_j)}
			{} & {}\leq C h^{p+1-r}\norm{u}_{H^{s}(\Gamma_j)},\qquad 0\leq r \leq p+1,
		\end{align}
		for integers $r,s$. Note that $C$ in \cite[Corollary 5.12]{BuffaActa} depends on $p$, but we choose $C$ as the maximal constant arising for any $p+1\leq s$.

		This immediately translates to the desired assertion on the boundary with $r = 0$ by continuity of the pullbacks and subadditivity of the $H^r$-norm.

		To extend this assertion to $r=-1/2,$ one can employ standard duality arguments as done, for example, by \cite[Theorem 4.1.32]{sauter_2010aa}.
	\end{proof}
\end{lemma}
Two remarks are due. First of all, one should note that the classical convergence result of order $3/2$ for boundary element methods corresponding to \eqref{problem::laplace} and \eqref{problem::helmholtz} for lowest order basis functions, i.e.~$p=0$, cf.~\cite{,Rjasanow_2007aa,sauter_2010aa,steinbach_2013aa}, manifests as a special case of the assertion above.

Moreover, by choosing $s$ maximally and choosing a sufficiently large $C$, for $p+1\leq s$ the term $\norm{u}_{\tilde H^{s}(\Gamma)}$ does not depend on $p$. Thus, in the case of degree elevation, one can expect exponential convergence up to the point where $p+1< s$, i.e., one is only limited by the regularity of the solution. For stronger assertions in dependence of $p,$ we refer to \cite{Veiga_2011aa}.

The result of Lemma \ref{cor::Approxcor} immediately transfers to the weak formulation \eqref{eq::varformdisc} and yield quasi-optimal convergence (in the $H^{-1/2}$ sense) of the unknown density.

However, for practical applications, often the pointwise error of the potential, i.e., the error of the solution in the interior/exterior of $\Omega$, is of greater relevance.
To this end, one can expect higher orders of convergence due the well known Aubin-Nitsche Lemma, as given, for example in \cite{Ciarlet_2002aa}.

\begin{lemma}[Aubin-Nitsche]\label{lem::aubinnitsche}
	For notational purposes, set $X \coloneqq H^{-p-2}(\Gamma)$. Let $w\in H^{-1/2}(\Gamma)$ and $w_h\in \S_{p,m}(\Gamma)$ be the solutions to \eqref{eq::varform} and \eqref{eq::varformdisc}, respectively. Then, one finds that
	\begin{align}
		\norm{w-w_h}_{X} \leq C \norm{w-w_h}_{H^{-1/2}(\Gamma)} \sup_{g\in X}\left(\norm{g}^{-1}_X\inf_{x_h\in \S_{p,m}(\Gamma)}\norm{w - x_h}_{H^{-1/2}(\Gamma)}\right).
	\end{align}
\end{lemma}

This leads immediately to the following assertion.

\begin{corollary}[Convergence of the potential]\label{cor::pwconv}
	Let $w$ and $w_h$ be as in Lemma~\ref{lem::aubinnitsche}. Moreover, assume that $w\in \tilde H^{p+1}(\Gamma)$
	and that $\Gamma$ is a boundary of class $C^{p+1,1}$ such that the space $H^{p+2}(\Gamma)$ is
	well defined.
	For all $\bb x\notin\Gamma$, there exists a positive constant $C=C(\bb x)$ 
	such that
	\begin{align}
		\abs{u(\bb x)-u_h(\bb x)}\leq  Ch^{2p+3}\norm{w}_{\tilde H^{p+1}(\Gamma)}\label{eq::pwconv}
	\end{align}
	holds for the functions $u$ and $u_h$ obtained by inserting $w$ and $w_h$ into \eqref{eq::repform}, respectively.
	\begin{proof}
		By the representation formula, the Aubin-Nitsche Lemma and estimate \eqref{eq::estimate} it holds
		\begin{align}
			{}     &{}\abs{u(\bb x)-u_h(\bb x)}\nonumber\\
			&\quad={}     {}\abs{\int_\Gamma u^*(\bb x,\bb y) (w(\bb y)-w_h(\bb y))\opd \bb y}\nonumber\\
			&\quad\leq{} {} \norm{u^*(\bb x,\cdot)}_{H^{p+2}(\Gamma)}\norm{w-w_h}_{H^{-p-2}(\Gamma)}\label{eq::smoothkernel}\\
			&\quad\leq{}  {}C \norm{w-w_h}_{H^{-1/2}(\Gamma)} \sup_{g\in  X}\left(\norm{g}^{-1}_{X}\inf_{v_h\in \S_{p,m}(\Gamma)}\norm{w - v_h}_{H^{-1/2}(\Gamma)}\right)\label{eq::constantcomesin} \\
			&\quad\leq {}  {}C h^{2p+3}\norm{w}_{\tilde H^{p+1}(\Gamma)},\nonumber
		\end{align}
		since $u^*(\bb x,\cdot)\in H^{s}(\Gamma)$ for arbitrary $s$ and $\bb x\notin\Gamma$.
	\end{proof}
\end{corollary}

\begin{remark}\label{rem::nastyparam}
	Note that, for general Lipschitz polyhedra, functions of the regularity required in \eqref{eq::smoothkernel}
	might not even be definable. Moreover, also for smooth geometries, the regularity assumption might be violated due to irregular geometry mappings, cf.~Remark \ref{rem::smoothkernel}.

	Nonetheless, even for nonsmooth geometries, an increase in the rate of convergence similar to the corollary above can be observed, as long as the patchwise mappings are regular, as seen in Section \ref{sec::numex}.
\end{remark}

\FloatBarrier

\section{Implementation} \label{sec::impl}

The code used is based on the implementation as proposed in \cite{doelz2016} and has been enriched by the following features:
\begin{enumerate}
	\item NURBS mappings for the geometry description, including Bézier extraction \cite{Borden_2011aa}, a method for fast NURBS evaluation by considering B-splines locally as given by Bernstein polynomials, similar to the technique of Section \ref{subsec::embedding}. The algorithms have been vectorised via OpenMP \cite{openmp08}.
	\item A more general superspace projection (cf.~Section \ref{subsec::embedding}) for higher order polynomial ansatz functions is now used to generate the B-spline spaces as in \eqref{def::space}.
  \item The implementation is now capable of solving the Helmholtz problem via a GMRES scheme, while the Laplace problem employs a CG method.
  \item Naturally, higher order quadrature formulae must be used for higher order methods.
  Since the order of quadrature must also be increased logarithmi\-cally towards the singularity arising from \eqref{eq::varformdisc}, Gauß quadrature has been implemented up to order 40. This allows the investigation of the effects of higher order methods, as will be explained in Section \ref{subsec::fmm}. Our approach is based on the one analysed in \cite{sauter_2010aa},
which has been adapted to our setting in \cite{harbrech2001}.
\end{enumerate}

To cope with the large, densely populated matrices (cf.~Table \ref{tab::mat::sph}), compression techniques must be utilised \cite{greengard2009fast,greengard1987fast,Kurz_2007aa}.
Algorithms, analysis and error estimates for the concepts to be reviewed in Sections \ref{subsec::clustertree}, \ref{subsec::fmm} and \ref{subsec::embedding} can be found in \cite{doelz2016}.

\subsection{Cluster Tree Structure and Localised Kernel Functions} \label{subsec::clustertree}

Each mapping of \eqref{def::geom} induces a patch $\Gamma_j$, to which uniform refinement (in terms of the reference domain) is applied. This procedure generates a nested sequence of meshes, which, for each level of refinement $m$, consists of $4^m$ elements per patch, on which the standard tensor product B-spline spaces can be defined, as in \eqref{def::tpspline}.

\begin{figure}
	\centering
	\begin{subfigure}{\textwidth}\centering
		\begin{tikzpicture}
			\draw (0,0) -- (0,2);
			\draw (0,2) -- (2,2);
			\draw (2,2) -- (2,0);
			\draw (2,0) -- (0,0);
			\node (C1) at (1,1) {$\Gamma_{i,0,0}$};

			\node (A) at (2.5,1) {};
			\node (B) at (3.5,1) {};
			\draw [->] (A.center) -- (B.center) node [midway, above] {refine};

			\draw (0+4,0) -- (0+4,2);
			\draw (0+4,2) -- (2+4,2);
			\draw (2+4,2) -- (2+4,0);
			\draw (2+4,0) -- (0+4,0);
			\draw (0+4,1) -- (2+4,1);
			\draw (0+4+1,0) -- (0+4+1,2);
			\node (C21) at (1+4-.5,1-.5) {$\Gamma_{i,1,0}$};
			\node (C22) at (1+4+.5,1-.5) {$\Gamma_{i,1,1}$};
			\node (C23) at (1+4+.5,1+.5) {$\Gamma_{i,1,2}$};
			\node (C24) at (1+4-.5,1+.5) {$\Gamma_{i,1,3}$};

			\node (C) at (6.5,1) {};
			\node (D) at (7.5,1) {};
			\draw [->] (C.center) -- (D.center) node [midway, above] {refine};

			\path [pattern = north east lines] (0+9,0) rectangle (0+9+1,1);

			\draw (0+8,0) -- (0+8,2);
			\draw (0+8,2) -- (2+8,2);
			\draw (2+8,2) -- (2+8,0);
			\draw (2+8,0) -- (0+8,0);
			\draw (0+8,1) -- (2+8,1);
			\draw (0+8+1,0) -- (0+8+1,2);
			\draw (0+8+.5,0) -- (0+8+.5,2);
			\draw (0+8+1.5,0) -- (0+8+1.5,2);
			\draw (0+8,1.5) -- (2+8,1.5);
			\draw (0+8,.5) -- (2+8,.5);

			\draw[ultra thick] (0+9,0) -- (0+9,1);
			\draw[ultra thick] (0+9+0,0) -- (0+9+1,0);
			\draw[ultra thick] (0+9+1,1) -- (0+9,1);
			\draw[ultra thick] (0+9+1,1) -- (0+9+1,0);

		\end{tikzpicture}
		\caption{Refinement of the patch induced by the $i$-th mapping. Bold region corresponds to cluster $\Gamma_{\bb \lambda}$ with $\bb \lambda = (i,1,1)$.}
	\end{subfigure}\\[.5cm]
	\begin{subfigure}{\textwidth}\centering
		\begin{tikzpicture}[level distance=1cm,level 1/.style={sibling distance=3.1cm},level 2/.style={sibling distance=2cm}]
			\node {$\Gamma_{i,0,0}$}
			child {node {$\Gamma_{i,1,0}$}
					child {node {$\Gamma_{i,2,0},\dots,\Gamma_{i,2,3}$}}
				}
			child {node {$\Gamma_{i,1,1}$}
					child {node {$\Gamma_{i,2,4},\dots,\Gamma_{i,2,7}$}}
				}
			child {node {$\Gamma_{i,1,2}$}
					child {node {$\Gamma_{i,2,8},\dots,\Gamma_{i,2,11}$}}}
			child {node {$\Gamma_{i,1,3}$}
					child {node {$\Gamma_{i,2,12},\dots,\Gamma_{i,2,15}$}}
				};
		\end{tikzpicture}
		\caption{Quadtree structure induced by refinement of the $i$-th mapping.}
	\end{subfigure}
	\caption{Refinement, labelling, and quadtree structure.}\label{fig::refinement}
\end{figure}
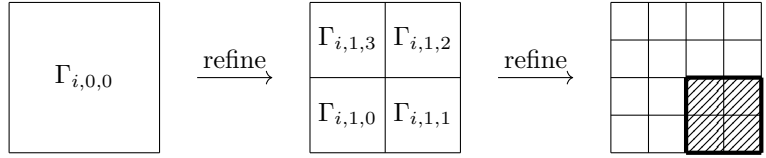
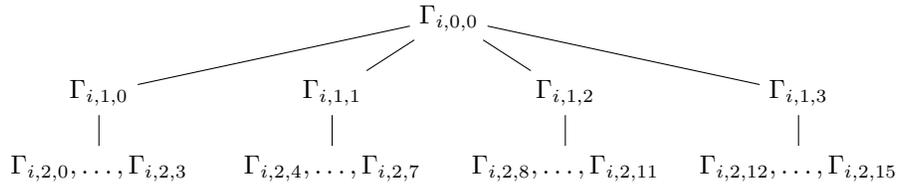

This mesh refinement strategy induces a quadtree structure on the geometry, cf.~Figure \ref{fig::refinement}. Each element $\Gamma_{i,j,k}$ within the nested sequence of meshes will be refered to by a tuple $(i,j,k)\eqqcolon \bb\lambda$, where $i$ denotes the corresponding parametric mapping, $j$ showcases the level of refinement of the element and $k$ denotes the index of the element in hierarchical order. For notational purposes, we will define $\abs{\bb \lambda}\coloneqq j.$
Each instance of $\Gamma_{\bb\lambda}$ is considered as a \emph{cluster} in terms of the fast multipole method, in the sense that $\Gamma_{\bb\lambda}$ will be considered as the set of tree leafs appended to the subtree with root $\Gamma_{\bb\lambda}.$ Na\"ively said, $\Gamma_{\bb \lambda}$ can be visualised as ``a square region on the geometry''.
The hierarchical ordered collection of all $\Gamma_{\bb\lambda}$ will be called \emph{cluster tree} and denoted by $\mathcal T$.

For each pair of clusters in $\mathcal T$, the fundamental solution of either \eqref{problem::laplace} or \eqref{problem::helmholtz} can be localised to a \emph{localised kernel function}
\begin{align}
	k_{\bb\lambda,\bb \lambda'} \colon [0,1]^2\times [0,1]^2 & \to\mathbb K, 
\end{align}
which corresponds to an evaluation of \eqref{eq::fundament1} or \eqref{eq::fundament2} at the points
$\bb x\in \Gamma_{\bb\lambda}$
and
$\bb y\in \Gamma_{\bb\lambda'},$
reparametrised to $[0,1]^2\times [0,1]^2$ according to the given geometry mapping and its surface measure.
This reduces the dimension (in terms of input variables) of the fundamental solution artificially.
Explicit construction of the localised kernel function is straightforward, but very technical, and can be found in \cite[Sec.~2]{doelz2016}.

\begin{remark}\label{rem::smoothkernel}
	This approach reduces the computational effort of the method considerably.
	However, due to the localisation of the kernel through the geometry mappings, the kernel function may lose regularity due to nonsmooth parametrisations, which might reinforce effects as discussed in Remark \ref{rem::nastyparam}.
	In the case of knot repetition, this could impede the quality of the interpolation, causing loss of convergence since a regularity assumption on the interpolated kernel as in equation \eqref{eq::smoothkernel} of Corollary \ref{cor::pwconv} requires global regularity.

	By reparametrisation, in which the $k$-th partial derivatives of the geometry mapping become nonsmooth, one can spoil the convergence rates achieved in potential evaluation. This effect is observable, although the parametrised geometry is still geometrically smooth, cf.~Section \ref{subsec::perturbed}.
\end{remark}

At this point, one should state that this problem is not specific to the method investigated, and smoothness of the geometry mappings in all higher order methods suffers from similar effects, if one were to remove regularity requirements in \eqref{def::geom}.

\subsection{Fast Multipole Method} \label{subsec::fmm}

For the compression and reduction of computational complexity, we utilise an interpolation-based multipole method.

Roughly, one starts by introducing a matrix structure via an admissibility criterion, that divides the interaction of sets of basis functions in the physical domain into different categories, based on their respective distance to each other.
Based on this, the integral operator of \eqref{problem::laplace} or \eqref{problem::helmholtz}, respectively, will be approximated by a polynomial interpolation, yielding a degenerate kernel approach.

This is largely analogous to classical multipole approaches \cite{greengard2009fast,greengard1987fast}. However, instead of implementing a fixed expansion via spherical harmonics, the kernel is represented by a polynomial series which is generated through interpolation.

We will make this brief explaination more precise. For an in depth review of the concepts, we refer to \cite{doelz2016} and \cite{Hackbusch_2002aa}.

\begin{definition}
	Given the block cluster tree structure and a fixed parameter $0<\mu <1$, define admissible blocks given by all $ \Gamma_{\bb\lambda}\times\Gamma_{\bb\lambda'}$ such that
	\begin{align}
		\max\lbrace \operatorname{diam}(\Gamma_{\bb\lambda}), \operatorname{diam}(\Gamma_{\bb\lambda'})\rbrace <\mu \operatorname{dist(\Gamma_{\bb\lambda},\Gamma_{\bb\lambda^ \prime})}.
	\end{align}
	The largest set of non-admissible blocks $\Gamma_{\bb \lambda}\times \Gamma_{\bb \lambda'}$ will be called the \emph{far field}, while the remaining non-admissible blocks will be called the \emph{near field}.
\end{definition}

Fixing tensor product Lagrange polynomials $L_{\bb m}(\hat {\bb x}) = L_{m_1}(x)L_{m_2}(y)$ for $0\leq m_1,m_2\leq p$ and corresponding interpolation points $\bb z_{\bb m} \coloneqq (z_{m_1},z_{m_2}),$
one can utilise an approximation of the singular kernels arising from \eqref{eq::varformdisc} via a degenerate kernel of the form
\begin{align}
	k_{\bb\lambda,\bb\lambda'}(\hat{\bb x},\hat{\bb y}) & \approx \sum_{\norm{\bb m}_{\infty},\norm{\bb m'}_{\infty}\leq p} k_{\bb\lambda,\bb\lambda'}(\bb z_{\bb m},\bb z_{\bb m'})L_{\bb m}(\hat{\bb x})L_{\bb m'}(\hat{\bb y}).
\end{align}
Setting $I=[0,1],$ the corresponding matrix block entries are then of the form
\begin{align}
    &[A_{\bb \lambda,\bb \lambda'}]_{i,j}\nonumber\\
   & \approx{} {} \iint_{I^2}\sum_{\norm{\bb m}_\infty,\norm{\bb m'}_\infty\leq p} k_{\bb\lambda,\bb\lambda'}(\hat{\bb z}_{\bb m},\hat{\bb z}_{\bb m'})L_{\bb m}(\hat{\bb x})\hat \phi_i(\hat{\bb x})L_{\bb m'}(\hat{\bb y})\hat \phi_j(\hat{\bb y})\opd \hat{\bb x}\opd \hat{\bb y}\\
  &= {} {}\sum_{\norm{\bb m}_{\infty},\norm{\bb m'}_\infty\leq p} k_{\bb\lambda,\bb\lambda'}(\hat{\bb z}_{\bb m},\hat{\bb z}_{\bb m'})\int_{I^2} L_{\bb m}(\hat{\bb x})\hat \phi_i(\hat{\bb x})\opd \hat{\bb x}\int_{I^2}  L_{\bb m'}(\hat{\bb y})\hat \phi_j(\hat{\bb y})\opd \hat{\bb y}.\nonumber
\end{align}
The integral terms are independent of the parametrisation, which appears only in the form of the interpolation values $k_{\bb\lambda,\bb\lambda'}(\hat{\bb x}_{\bb m},\hat{\bb x}_{\bb m'})$. Thus, apart from the kernel evaluation, the term above can be stored in moment matrices independent of the parametrisation. Utilising the tensor product structure of the mapping, the moment matrices can be simplified even further, cf.~\cite[Sec.~5.2]{doelz2016}.

Moreover, utilising the concept of {nested cluster bases} \cite{Hackbusch_2002aa}, which amounts to an $\mathcal H^2$-matrix representation, an explicit computation of the moment matrices can be avoided for almost all matrix blocks \cite[Sec.~5.4]{doelz2016}.

\begin{remark}
	It should be noted that, for this compression approach to work, one merely needs to require the kernel to be \emph{analytically standard}, a technical notion defined in the literature cited above. This requirement is much weaker compared to other requirements from different compression approaches, cf.~\cite{Hackbusch_2002aa,Kurz_2007aa}. We omit refering to this requirement, because it is fulfilled by both, \eqref{eq::fundament1} and \eqref{eq::fundament2}.
\end{remark}

\subsection{Approach via Discontinuous Superspace} \label{subsec::embedding}

Bézier extraction, as used for evaluation of the geometry, can be used to evaluate the basis functions as well. However, in a na\"ive implementation, this might require additional evaluations of the kernel and the geometry, which is a computational effort that can be omitted by constructing the matrix $\bb A$ corresponding to \eqref{eq::varformdisc} element wise. The functions in $\S_{p,m}(\Gamma)$ are interpreted as a linear combination of Bernstein polynomials defined on the individual elements, i.e.~the quadrilatiral domains given by the refinement process detailed above. This can be formalized as follows.

Let $n$ denote the number of patches representing $\Gamma$ and fix a polynomial degree $p$.
For $0\leq i\leq p$, defining the Bernstein polynomials as
\begin{align}
	B_{i,p}(x)\coloneqq \dbinom{p}{i}x^i(1-x)^{p-i}, \qquad0\leq x\leq 1,
\end{align}
we apply on each element of the given refinement level $m$, w.l.o.g.~given by $[0,1]^2$, a tensor product structure as in
\begin{align}
	B_{i,j,p}(x,y) = B_{i,p}(x)\cdot B_{j,p}(y),\qquad & (x,y)\in [0,1].
\end{align}
By affine transformation and utilisation of the pullbacks induced by the $\bb F_j$, this yields a global discrete function space $\S^*_{p,m}(\Gamma)$ of dimension $k\coloneqq n\cdot((2^m)\cdot(p+1))^2,$ where $m$ again denotes the level of refinement\footnote{$2^m$ elements, with $(p+1)^2$ local Bernstein polynomials per element, on $n$ patches.}.
Since B-splines are piecewise polynomials, it clearly holds that $\S_{p,m}(\Gamma)\subseteq\S^*_{p,m}(\Gamma)$.

Assembly of the system matrix $\bb A$ corresponding to \eqref{eq::varform} w.r.t.~$\S^*_{p,m}(\Gamma)$ yields a large matrix of size $k\times k.$
However, due to the highly local support of the ansatz functions in $\S^*_{p,m}(\Gamma)$, the matrix $\bb A$ is handled easily by the multipole method reviewed in Section \ref{subsec::fmm}.

Having $\bb A$ in compressed format, we then apply the method by utilisation of a projection $\S_{p,m}(\Gamma)\to \S^*_{p,m}(\Gamma)$.
Said operator can be instantiated by application of a sparse transformation matrix $\bb T\in \R^{\ell\times k},$ with $\ell = n\cdot((2^m+p)^2) $ denoting the dimension of $\S_{p,m}(\Gamma)$. One thus finds
\begin{align}
	\bb T \bb A \bb T^\top \bb w = \bb T \bb g^*,
\end{align}
$\bb w$ hereby denoting the coefficients for the density in $\S_{p,m}(\Gamma)$ and $\bb g^*$ the coefficient vector corresponding to the right hand side $g$ w.r.t.~$\S^*_{p,m}(\Gamma)$.
Depending on $h$ and $p$, the corresponding values for $k$ can be found in Table \ref{tab::mat::sph}.

To summarise, the integral operator is discretised as a matrix w.r.t.\ $\S^*_{p,m}(\Gamma)$, to which the interpolation-based multipole method is applied.
Only afterwards, the relation between $\S^*_{p,m}(\Gamma)$ and $\S_{p,m}(\Gamma)$ is taken into account, which corresponds to application of a geometry independent sparse matrix and thus can be handled efficiently.

This approach offers huge flexibility in the implementation, for a reasonable price, since the utilisation of the compression techniques discussed in Section \ref{subsec::fmm} reduces the computational effort and the memory consumption considerably \cite[Sec.~5.3]{doelz2016}.

\begin{table}[!ht]
	\centering
	\caption{Dimension of the space $\S^*_{p,m}(\Gamma)$.}
	\label{tab::mat::sph}
	\begin{tabular}{l|l|l|l|l|l}\multicolumn{5}{c}{}\\
		      & \multicolumn{5}{c}{\textbf{Sphere}}\\\hline
		      & $m = 1$ & $m = 2$ & $m = 3$ & $m = 4$ & $m = 5$ \\\hline
		$p=0$ & 24      & 96      & 384     & 1,536   & 6,144   \\\hline
		$p=1$ & 96      & 384     & 1,536   & 6,144   & 24,576  \\\hline
		$p=2$ & 216     & 864     & 3,456   & 13,824  & 55,296  \\\hline
		$p=3$ & 384     & 1,536   & 6,144   & 24,576  & 98,304  \\\hline
		$p=4$ & 600     & 2,400   & 9,600   & 38,400  & 153,600 \\\hline\\
		      & \multicolumn{5}{c}{\textbf{Torus}}\\\hline
		      & $m = 1$ & $m = 2$ & $m = 3$ & $m = 4$ & $m = 5$ \\\hline
		$p=0$ & 64      & 256     & 1,024   & 4,096   & 16,384  \\\hline
		$p=1$ & 256     & 1,024   & 4,096   & 16,348  & 65,536  \\\hline
		$p=2$ & 576     & 2,304   & 9,216   & 36,864  & 147,456 \\\hline
		$p=3$ & 1,024   & 4,096   & 16,384  & 65,536  & 262,144 \\\hline
		$p=4$ & 1,600   & 6,400   & 25,600  & 102,400 & 409,600 \\\hline\\
		      & \multicolumn{5}{c}{\textbf{Fichera cube}}\\\hline
		      & $m = 1$ & $m = 2$ & $m = 3$ & $m = 4$ & $m = 5$ \\\hline
		$p=0$ & 96      & 384     & 1,536   & 6,144   & 24,576  \\\hline
		$p=1$ & 384     & 1,536   & 6,144   & 24,576  & 98,304  \\\hline
		$p=2$ & 864     & 3,456   & 13,824  & 55,296  & 221,184 \\\hline
		$p=3$ & 1,536   & 6,144   & 24,576  & 98,304  & 393,216 \\\hline
		$p=4$ & 2,400   & 9,600   & 38,400  & 153,600 & 614,400 \\\hline
	\end{tabular}
\end{table}

\section{Numerical Examples} \label{sec::numex}

\begin{figure}[h]
	\centering
	\begin{subfigure}{.31\textwidth}
		\includegraphics[width=\textwidth]{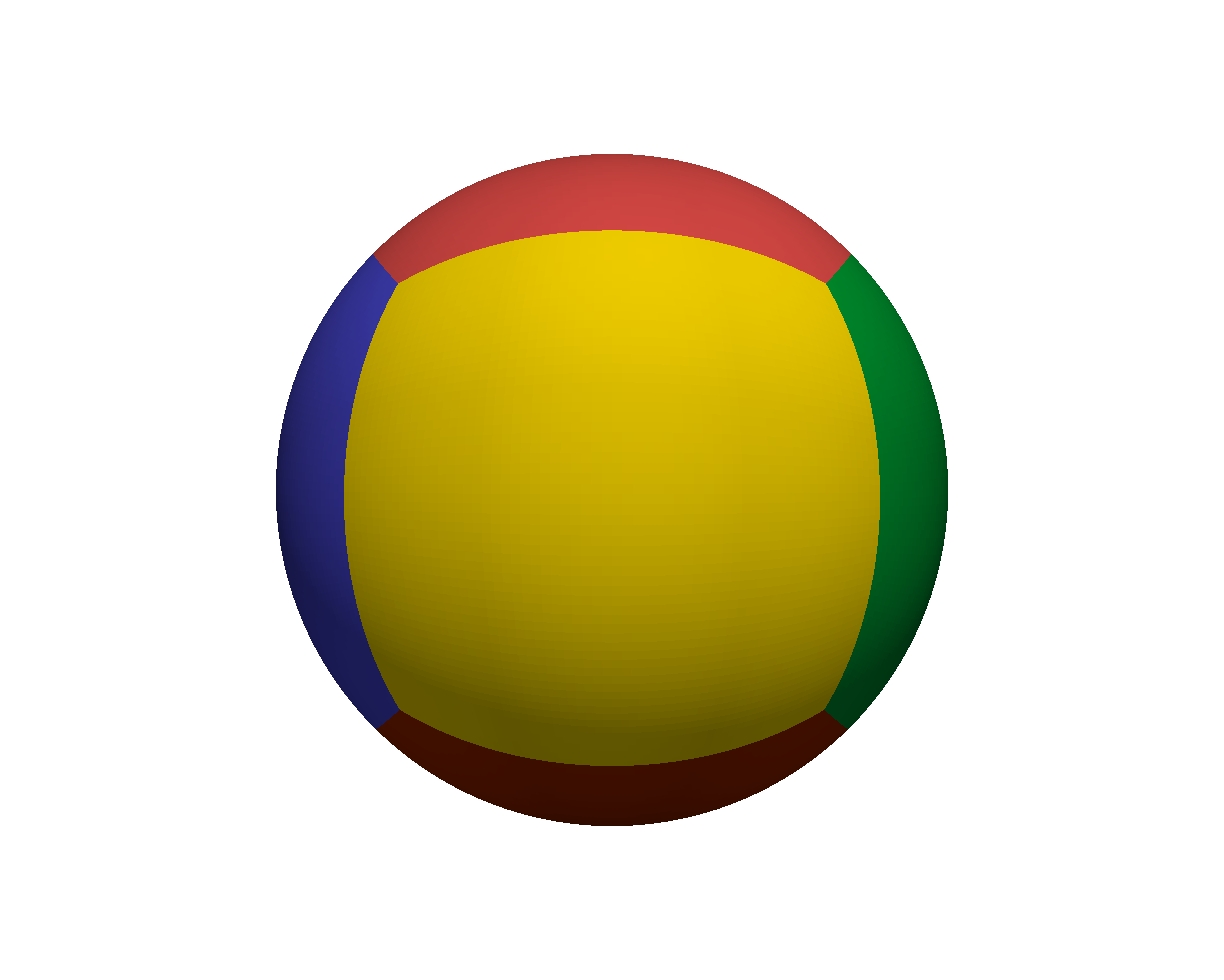}
		\label{fig::geom::sphere}
	\end{subfigure}
	\begin{subfigure}{.31\textwidth}
		\centering
		\includegraphics[width=\textwidth]{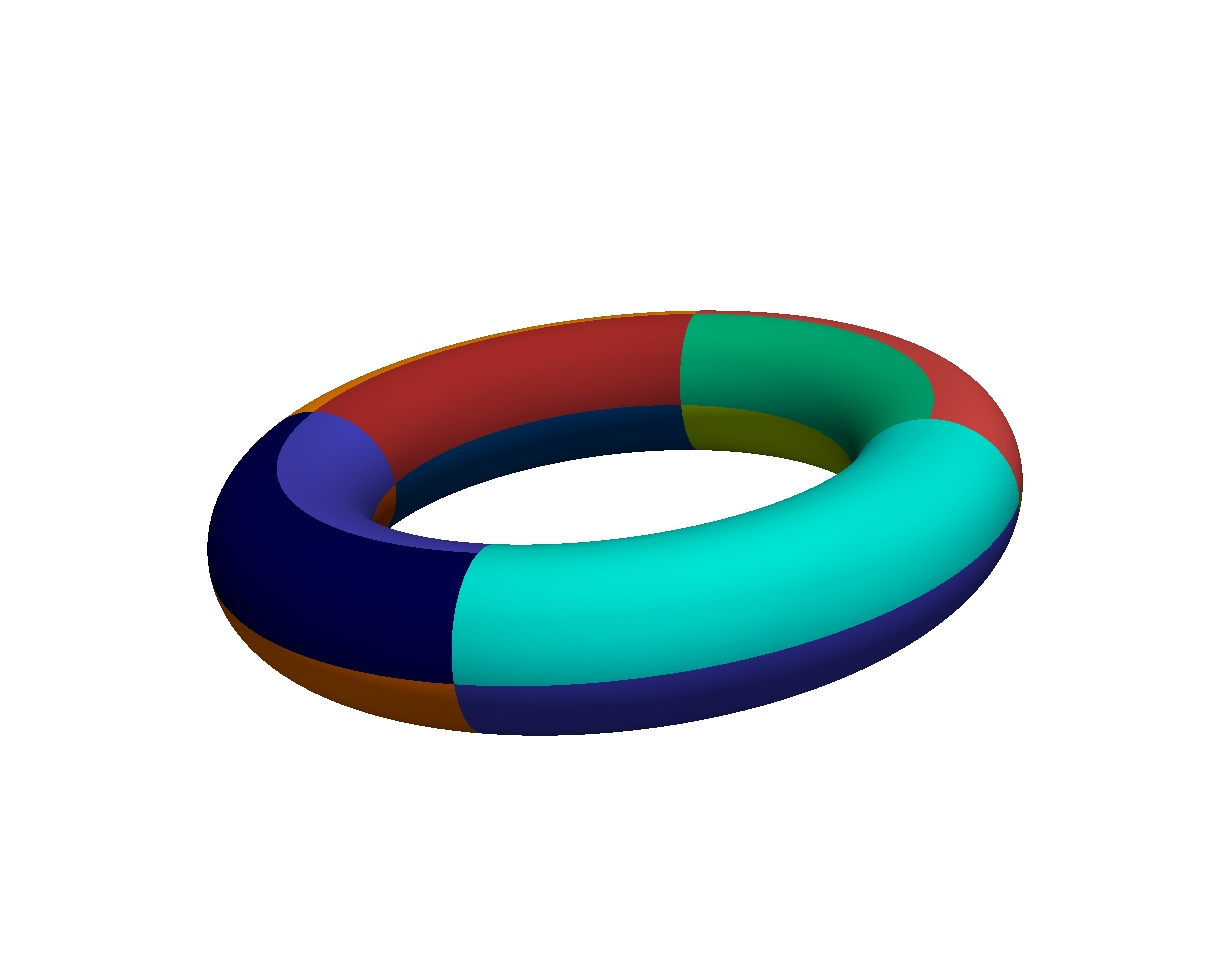}
		\label{fig::geom::torus}
	\end{subfigure}
	\begin{subfigure}{.31\textwidth}
		\centering
		\includegraphics[width=\textwidth]{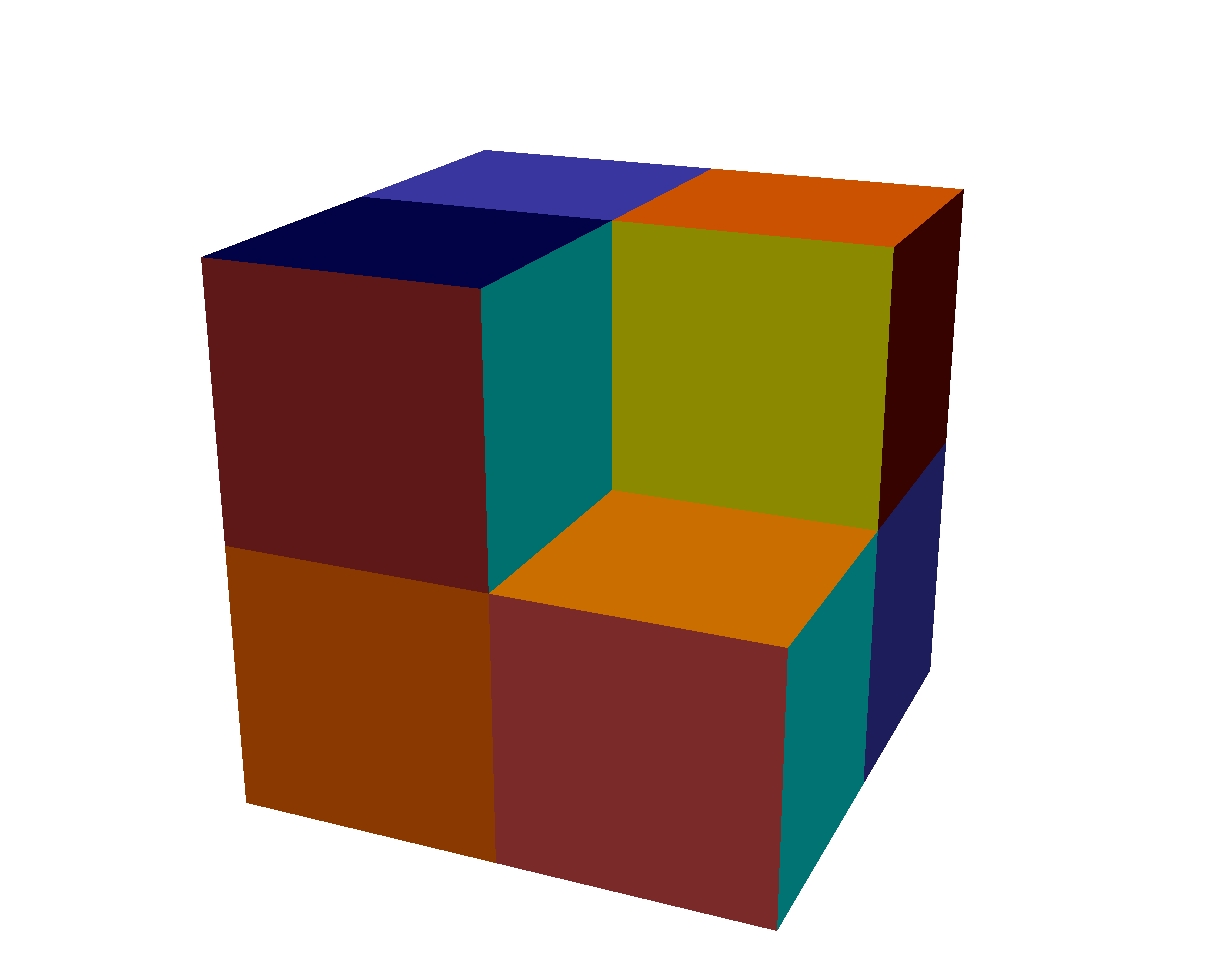}
		\label{fig::geom::fichera}
	\end{subfigure}
	\caption{Parametrisation of the geometries. Sphere, torus, and Fichera cube.}
	\label{fig::geom}
\end{figure}

As test geometries we chose to utilise basic geometries for simple comparison with different methods, cf.~Figure \ref{fig::geom}.
The sphere geometry is the unit sphere and has been parametrised using 6 NURBS patches of degree $p=4$ arranged in a cube topology.
For the torus geometry, 16 patches of degree $p=2$ are utilised. The outer radius of the torus is 2, while the inner radius is 0.5. Each quarter of the torus is given as four patches, each forming a pipe performing a  $90^{\circ}$ bent.

Lastly, the Fichera cube is given as a geometry consisting of 24 linear patches arranged as a cube of length 2, with an indented corner.

All potential errors are to be understood as total, pointwise errors in complete analogy to the estimate of Corollary \ref{cor::pwconv}.

\subsection{Laplace Problem}\label{subsec::laplnumex}

Numerical validation has been done through the method of manufactured solutions, i.e.~by taking boundary data given by a known function satisfying all requirements to be a valid solution to the problem investigated. In the Laplace case, the Dirichlet boundary data was generated by the spherical harmonic
\begin{align}
	Y^2_0(x,y,z) = \sqrt{\frac{5}{16\pi}}(3z^2-1).
\end{align}
On the sphere, the analytical solution of the density is known \cite[Sec.~7.3]{doelz2016}, and the $L^2$-error of the density has been computed as well, which is known to converge in accordance to the rates predicted in Corollary \ref{cor::Approxcor}. The evaluation points of the potential within the interior have been selected as follows: The interior of the domain is covered by cubes of lenght 0.05, and cubes with a distance of less than 0.15 to the boundary are discarded, cf.~Figure \ref{fig::spherepic}. Then, on all corners of the cubes, the potential is evaluated and the maximum error is chosen.

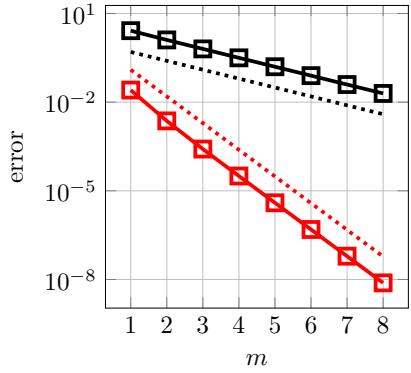
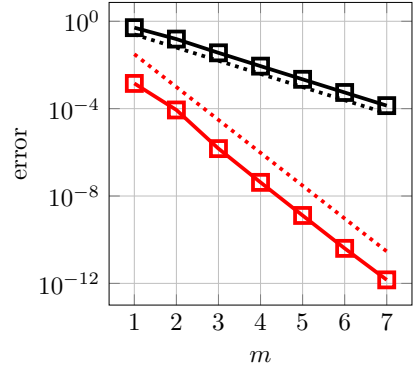
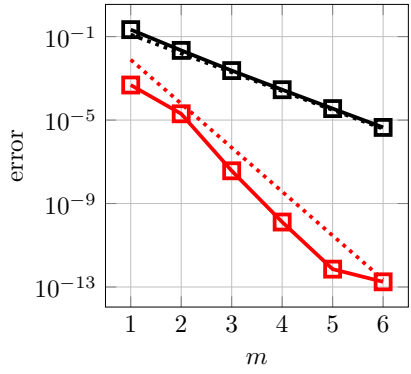
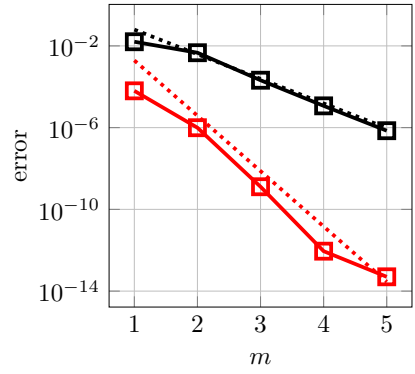
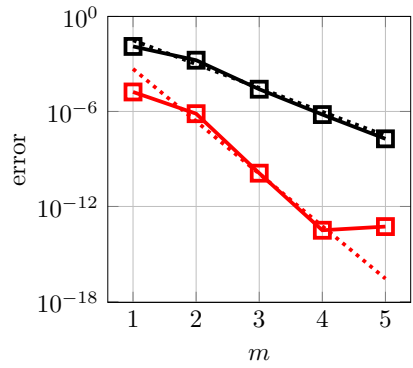
\begin{figure}
	\centering
	\begin{subfigure}{.48\textwidth}
		\begin{tikzpicture}[scale = .95]
			\begin{semilogyaxis}[
					width=1\columnwidth,
					height=1\columnwidth,
					xlabel={$m$},
					ylabel={error},
					xticklabels={1,2,3,4,5,6,7,8},
					xtick = {1,2,3,4,5,6,7,8},
					grid=major,
				]
				\addplot[line width=1.5pt,mark = none,red,dotted]  table[x index=0, y index=2] {data/lapl/sphere/sl0b_sphere_1_6_12};
				\addplot[line width=1.5pt,mark = none,dotted]  table[x index=0, y index=4] {data/lapl/sphere/sl0b_sphere_1_6_12};
				\addplot[line width=1.5pt,mark = square,mark options={solid},mark size=3,red]  table[x index=0, y index=3] {data/lapl/sphere/sl0b_sphere_1_6_12};
				\addplot[line width=1.5pt,mark = square,mark options={solid},mark size=3] table[x index=0, y index=5] {data/lapl/sphere/sl0b_sphere_1_6_12};
			\end{semilogyaxis}
		\end{tikzpicture}
		\caption{$p=0$}
	\end{subfigure}
	\begin{subfigure}{.48\textwidth}
		\begin{tikzpicture}[scale = .95]
			\begin{semilogyaxis}[
					width=1\columnwidth,
					height=1\columnwidth,
					xlabel={$m$},
					ylabel={error},
					xticklabels={1,2,3,4,5,6,7},
					xtick = {1,2,3,4,5,6,7},
					grid=major,
				]
				\addplot[line width=1.5pt,mark = square,mark options={solid},mark size=3,red] table[x index=0, y index=3] {data/lapl/sphere/sl1b_sphere_1_6_12};
				\addplot[line width=1.5pt,mark = none,red,dotted]  table[x index=0, y index=2] {data/lapl/sphere/sl1b_sphere_1_6_12};
				\addplot[line width=1.5pt,mark = none,dotted]  table[x index=0, y index=4] {data/lapl/sphere/sl1b_sphere_1_6_12};
				\addplot[line width=1.5pt,mark = square,mark options={solid},mark size=3] table[x index=0, y index=5] {data/lapl/sphere/sl1b_sphere_1_6_12};
			\end{semilogyaxis}
		\end{tikzpicture}
		\caption{$p=1$}
	\end{subfigure}\\[.5cm]
	\begin{subfigure}{.48\textwidth}
		\begin{tikzpicture}[scale = .95]
			\begin{semilogyaxis}[
					width=1\columnwidth,
					height=1\columnwidth,
					xlabel={$m$},
					ylabel={error},
					xticklabels={1,2,3,4,5,6},
					xtick = {1,2,3,4,5,6},
					grid=major,
				]
				\addplot[line width=1.5pt,mark = square,mark options={solid},mark size=3,red] table[x index=0, y index=3] {data/lapl/sphere/sl2b_sphere_1_6_12};
				\addplot[line width=1.5pt,mark = none,red,dotted]  table[x index=0, y index=2] {data/lapl/sphere/sl2b_sphere_1_6_12};
				\addplot[line width=1.5pt,mark = none,dotted]  table[x index=0, y index=4] {data/lapl/sphere/sl2b_sphere_1_6_12};
				\addplot[line width=1.5pt,mark = square,mark options={solid},mark size=3] table[x index=0, y index=5] {data/lapl/sphere/sl2b_sphere_1_6_12};
			\end{semilogyaxis}
		\end{tikzpicture}
		\caption{$p=2$}
	\end{subfigure}
	\begin{subfigure}{.48\textwidth}
		\begin{tikzpicture}[scale = .95]
			\begin{semilogyaxis}[
					width=1\columnwidth,
					height=1\columnwidth,
					xlabel={$m$},
					ylabel={error},
					xticklabels={1,2,3,4,5,6,7},
					xtick = {1,2,3,4,5,6,7},
					grid=major,
				]
				\addplot[line width=1.5pt,mark = square,mark options={solid},mark size=3,red] table[x index=0, y index=3] {data/lapl/sphere/sl3b_sphere_1_6_12};
				\addplot[line width=1.5pt,mark = none,red,dotted]  table[x index=0, y index=2] {data/lapl/sphere/sl3b_sphere_1_6_12};
				\addplot[line width=1.5pt,mark = none,dotted]  table[x index=0, y index=4] {data/lapl/sphere/sl3b_sphere_1_6_12};
				\addplot[line width=1.5pt,mark = square,mark options={solid},mark size=3] table[x index=0, y index=5] {data/lapl/sphere/sl3b_sphere_1_6_12};
			\end{semilogyaxis}
		\end{tikzpicture}
		\caption{$p=3$}
	\end{subfigure}\\[.5cm]
	\begin{subfigure}{.48\textwidth}
		\begin{tikzpicture}[scale = .95]
			\begin{semilogyaxis}[
					width=1\columnwidth,
					height=.96\columnwidth,
					xlabel={$m$},
					ylabel={error},
					xticklabels={1,2,3,4,5,6},
					xtick = {1,2,3,4,5,6},
					grid=major,
				]
				\addplot[line width=1.5pt,mark = square,mark options={solid},mark size=3,red] table[x index=0, y index=3] {data/lapl/sphere/sl4b_sphere_1_6_12};
				\addplot[line width=1.5pt,mark = none,red,dotted]  table[x index=0, y index=2] {data/lapl/sphere/sl4b_sphere_1_6_12};
				\addplot[line width=1.5pt,mark = none,dotted]  table[x index=0, y index=4] {data/lapl/sphere/sl4b_sphere_1_6_12};
				\addplot[line width=1.5pt,mark = square,mark options={solid},mark size=3] table[x index=0, y index=5] {data/lapl/sphere/sl4b_sphere_1_6_12};
			\end{semilogyaxis}
		\end{tikzpicture}
		\caption{$p=4$}
	\end{subfigure}

	\caption{Error for the interior Laplace problem of the sphere. Red: Maximal potential error. Black: $L^2(\Gamma)$ error of the density.  The red dotted line represents the ideal convergence rate of $h^{2p+3},$ in accordance to Corollary \ref{cor::pwconv}. The black dotted line shows the expected convergence of the $L^2$-error of the density of order $p$.}
	\label{fig::numres::lapl::sphere}
\end{figure}

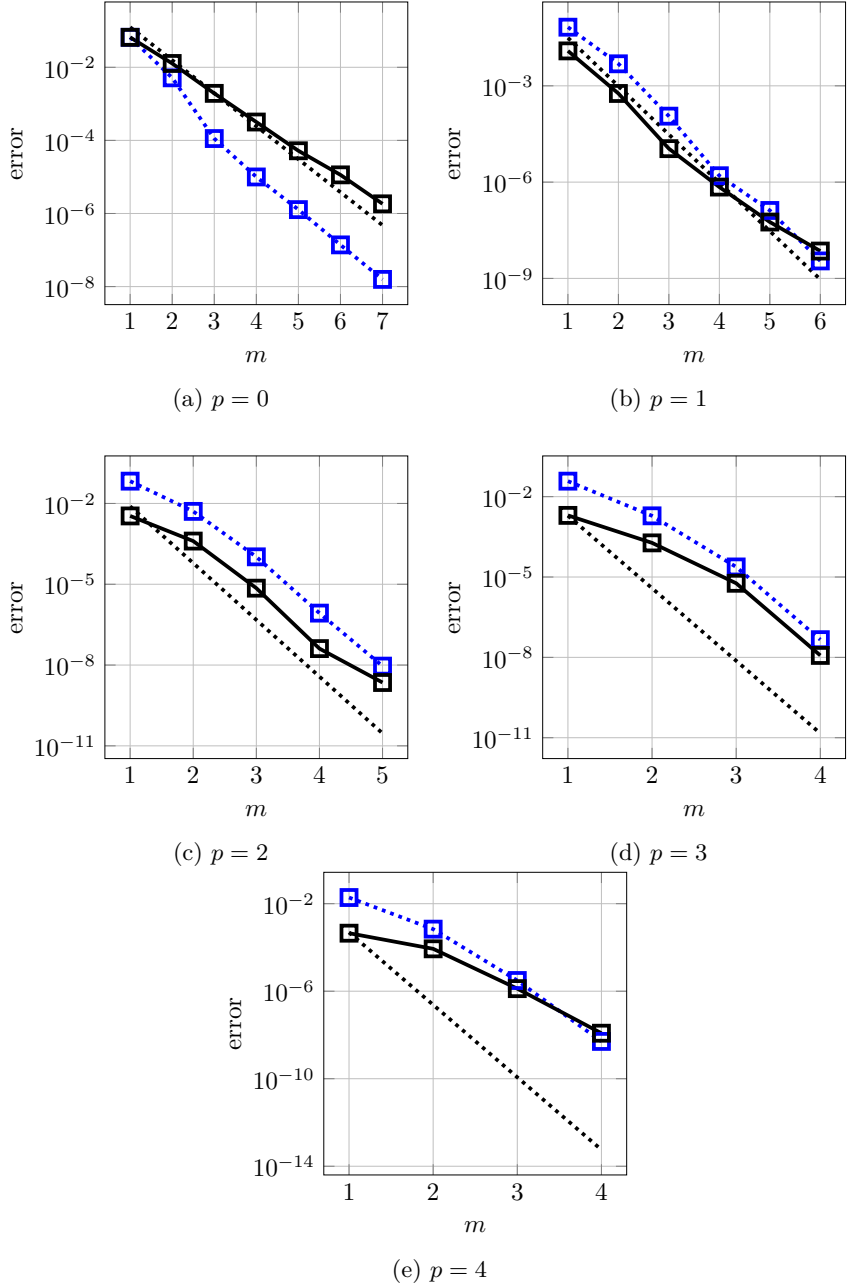
\begin{figure}
	\centering
	\begin{subfigure}{.48\textwidth}
		\begin{tikzpicture}[scale = .95]
			\begin{semilogyaxis}[
					width=1\columnwidth,
					height=1\columnwidth,
					xlabel={$m$},
					ylabel={error},
					xticklabels={1,2,3,4,5,6,7},
					xtick = {1,2,3,4,5,6,7},
					grid=major,
				]
				\addplot[line width=1.5pt,mark = square,mark options={solid},mark size=3,blue,dotted] table[x index=0, y index=3] {data/lapl/torus/sl0b_1_5_10};
				\addplot[line width=1.5pt,black,dotted] table[x index=0, y index=2] {data/lapl/torus/sl0b_1_5_10};
				\addplot[line width=1.5pt,mark = square,mark options={solid},mark size=3] table[x index=0, y index=3] {data/lapl/fichera/sl0b_fichera_1_5_10};
			\end{semilogyaxis}
		\end{tikzpicture}
		\caption{$p=0$}
	\end{subfigure}
	\begin{subfigure}{.48\textwidth}
		\begin{tikzpicture}[scale = .95]
			\begin{semilogyaxis}[
					width=1\columnwidth,
					height=1\columnwidth,
					xlabel={$m$},
					ylabel={error},
					xticklabels={1,2,3,4,5,6},
					xtick = {1,2,3,4,5,6},
					grid=major,
				]
				\addplot[line width=1.5pt,mark = square,mark options={solid},mark size=3,blue,dotted] table[x index=0, y index=3] {data/lapl/torus/sl1b_1_5_12};
				\addplot[line width=1.5pt,black,dotted] table[x index=0, y index=2] {data/lapl/torus/sl1b_1_5_12};
				\addplot[line width=1.5pt,mark = square,mark options={solid},mark size=3] table[x index=0, y index=3] {data/lapl/fichera/sl1b_fichera_1_5_12};
			\end{semilogyaxis}
		\end{tikzpicture}
		\caption{$p=1$}
	\end{subfigure}\\[.5cm]
	\begin{subfigure}{.48\textwidth}
		\begin{tikzpicture}[scale = .95]
			\begin{semilogyaxis}[
					width=1\columnwidth,
					height=1\columnwidth,
					xlabel={$m$},
					ylabel={error},
					xticklabels={1,2,3,4,5},
					xtick = {1,2,3,4,5},
					grid=major,
				]
				\addplot[line width=1.5pt,mark = square,mark options={solid},mark size=3,blue,dotted] table[x index=0, y index=3] {data/lapl/torus/sl2b_1_5_12};
				\addplot[line width=1.5pt,black,dotted] table[x index=0, y index=2] {data/lapl/torus/sl2b_1_5_12};
				\addplot[line width=1.5pt,mark = square,mark options={solid},mark size=3] table[x index=0, y index=3] {data/lapl/fichera/sl2b_fichera_1_5_12};
			\end{semilogyaxis}
		\end{tikzpicture}
		\caption{$p=2$}
	\end{subfigure}
	\begin{subfigure}{.48\textwidth}
		\begin{tikzpicture}[scale = .95]
			\begin{semilogyaxis}[
					width=1\columnwidth,
					height=1\columnwidth,
					xlabel={$m$},
					ylabel={error},
					xticklabels={1,2,3,4},
					xtick = {1,2,3,4},
					grid=major,
				]
				\addplot[line width=1.5pt,mark = square,mark options={solid},mark size=3,blue,dotted] table[x index=0, y index=3] {data/lapl/torus/sl3b_1_4_12};
				\addplot[line width=1.5pt,black,dotted] table[x index=0, y index=2] {data/lapl/torus/sl3b_1_4_12};
				\addplot[line width=1.5pt,mark = square,mark options={solid},mark size=3] table[x index=0, y index=3] {data/lapl/fichera/sl3b_ficher_1_4_12};
			\end{semilogyaxis}
		\end{tikzpicture}
		\caption{$p=3$}
	\end{subfigure}
	\begin{subfigure}{.48\textwidth}
		\begin{tikzpicture}[scale = .95]
			\begin{semilogyaxis}[
					width=1\columnwidth,
					height=1\columnwidth,
					xlabel={$m$},
					ylabel={error},
					xticklabels={1,2,3,4},
					xtick = {1,2,3,4},
					grid=major,
				]
				\addplot[line width=1.5pt,mark = square,mark options={solid},mark size=3,blue,dotted] table[x index=0, y index=3] {data/lapl/torus/sl4b_1_4_12};
				\addplot[line width=1.5pt,black,dotted] table[x index=0, y index=2] {data/lapl/torus/sl4b_1_4_12};
				\addplot[line width=1.5pt,mark = square,mark options={solid},mark size=3] table[x index=0, y index=3] {data/lapl/fichera/sl4b_fichera_1_4_12};
			\end{semilogyaxis}
		\end{tikzpicture}
		\caption{$p=4$}
	\end{subfigure}
	\caption{Potential error for the interior Laplace problem. Torus (blue, dotted) and Fichera cube (black). Maximal potential error. The dotted line represents the ideal convergence rate of $h^{2p+3},$ in accordance to Corollary \ref{cor::pwconv}.}
	\label{fig::lapl::fichera}
\end{figure}

The convergence plots in Figures \ref{fig::numres::lapl::sphere} and \ref{fig::lapl::fichera} confirm the expected convergence behaviour predicted by Corollaries \ref{cor::Approxcor} and \ref{cor::pwconv} for smooth geometries, where for larger $p$ we soon reach machine accuracy for the potential error of the sphere. 


One can also see that for large $m$ the convergence rate of the Fichera cube starts to fluctuate in the case of $p=1,2$. For $p=4$, it does not reach the ideal order of convergence at all. This can be attributed to the missing boundary regularity required by Corollary \ref{cor::pwconv}. However, one still sees a clear benefit from an increase of $p$.

Indeed, as can be seen in Table \ref{tab::time}, higher order approaches of $p\geq2$ prove advantagous in terms of time to solution over lower order approaches, in this case $p=0,1$. This can be attributed to the time saved in matrix assembly, since, for a set accuracy bound, the dimension of the space $\S_{p,m}(\Gamma)$ required to reach said bound is of far lower dimension for larger $p$.

\subsection{Perturbed Mappings}\label{subsec::perturbed}

To showcase the effect of perturbed mappings, the (patchwise) smooth mappings $\bb F_j\colon[0,1]^2\to \Gamma_j$ of the sphere geometry were perturbed componentwise by the function $f\colon[0,1] \to [0,1]$ given by 
\begin{align}
	f(x)\coloneqq\begin{cases}
		0.5\cdot x,&\text{ for } 0\leq x \leq 0.3,\\
		1.75 \cdot x - 0.375,&\text{ for } 0.3< x < 0.7,\\
		0.5 \cdot(x + 1),&\text{ for } 0.7\leq x \leq 1.
	\end{cases}
\end{align}
Note that $f\in H^1([0,1])$ but $f\notin H^2([0,1]).$ Application of such a perturbation will in general affect the regularity of the discrete space and the interpolation behaviour of the multipole method.

\begin{figure}
\centering
		\begin{tikzpicture}[scale = .95]
			\begin{semilogyaxis}[
					width=.7\columnwidth,
					height=.48\textwidth,
					xlabel={$m$},
					ylabel={error},
					xticklabels={1,2,3,4,5},
					xtick = {1,2,3,4,5},
					grid=major,
				]
				\addplot[line width=1.5pt,mark = triangle,mark options={solid},mark size=3,red] table[x index=0, y index=3] {data/pert/pert_a_o_1};
				
				\addplot[line width=1.5pt,mark = square,mark options={solid},mark size=3,red,dashed] table[x index=0, y index=3] {data/pert/pert_a_o_2};

				\addplot[line width=1.5pt,mark = diamond,mark options={solid},mark size=3,red,densely dotted] table[x index=0, y index=3] {data/pert/pert_a_o_3};
				\addplot[line width=1.5pt,black,dotted] table[x index=0, y index=2] {data/pert/pert_a_o_1};

				\legend{$p=0$,$p=1$,$p=2$};
			\end{semilogyaxis}
		\end{tikzpicture}
	\caption{Potential error for the interior Laplace problem with perturbed mappings and $p=0,1,2$. Maximal potential error. The dotted line represents a convergence rate of $h^3$ for the error.}
	\label{fig::lapl::pert}
\end{figure}
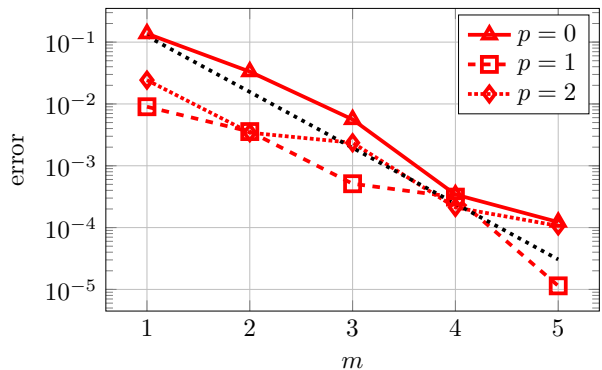

Figure \ref{fig::lapl::pert} shows the same numerical experiment as described in Section \ref{subsec::laplnumex} for the sphere, including a perturbation of the geometry mapping. Note that the geometry is still parametrised without an error, merely the regularity of the parametrisation is affected. One can see that, in this case, convergence rates for $p=1,2$ do not improve compared to the convergence rates achieved by a lowest order method.

\FloatBarrier

\subsection{Helmholtz Problem}\label{subsec::helmnumex}

For benchmarking the solution of Helmholtz problems within the exterior domain, the Dirichlet data has been induced by the function
\begin{align}
	f(\bb x)= \frac{\exp(i\kappa\norm{\bb{x}- \bb v})}{\norm{\bb x - \bb v}},
\end{align}
where $\bb v$ is a point within the domain ($(0,-2,0)$ for the torus, $(0.5,0.5,0.5)$ for all other geometries, and $\kappa = 0.5$)\footnote{The choices of $\kappa,$ i.e.~$0.5,$ $1,$ and $2$ correspond to audible frequencies in air, and are within the spectrum of practical relevance.}. For larger $\kappa$, note that one needs higher levels of refinement to discretise the resulting field. The function satisfies the Helmholtz equation and the Sommerfeld radiation condition in the exterior.

Potential evaluation in each case was done at $2^{m+2}$ points on a sphere with radius $6$ around the origin. They are distributed equally and they are nested, i.e.\ for an increase of $m$ new points are added to the set. 

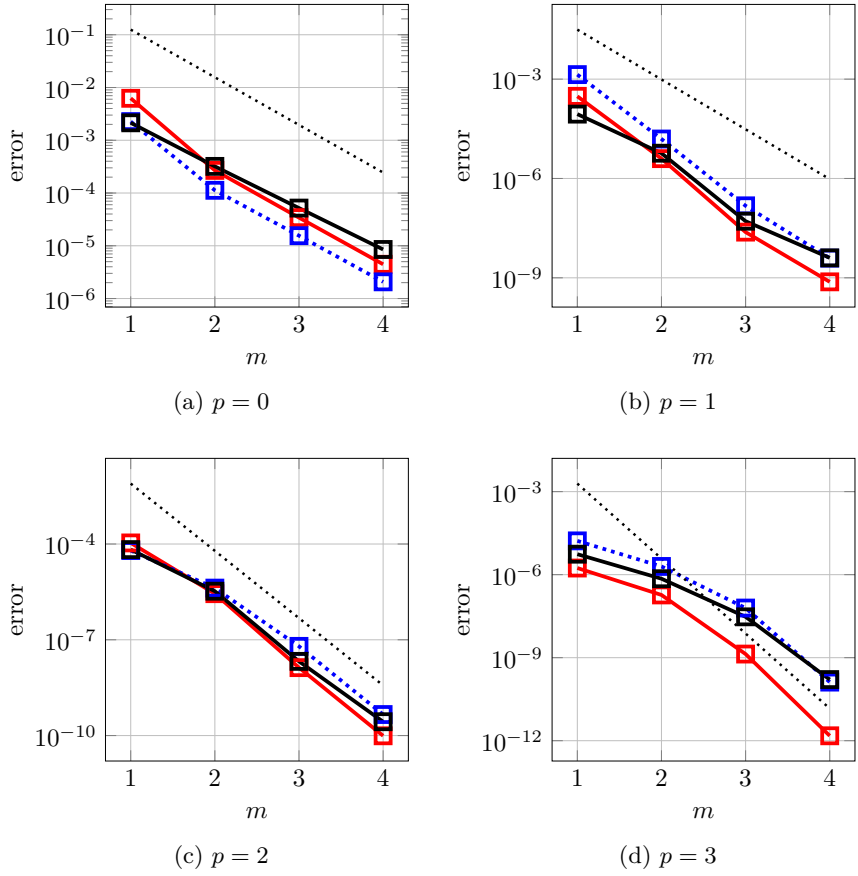
\begin{figure}[!h]
	\centering
	\begin{subfigure}{.48\textwidth}
		\begin{tikzpicture}[scale = .95]
			\begin{semilogyaxis}[
					width=1\columnwidth,
					height=1\columnwidth,
					xlabel={$m$},
					ylabel={error},
					xticklabels={1,2,3,4},
					xtick = {1,2,3,4},
					grid=major,
				]
				\addplot[line width=1.5pt,mark = square,mark options={solid},mark size=3,red]  table[x index=0, y index=3] {data/helm/tours/sh0b_t_1_4_12};
				\addplot[line width=1.5pt,mark = square,mark options={solid},dotted,mark size=3,blue]  table[x index=0, y index=3] {data/helm/sphere/sh0b_1_4_12};
				\addplot[line width=1.5pt,mark = square,mark options={solid},mark size=3]  table[x index=0, y index=3] {data/helm/ficher_data/sh0b_f_1_4_10};
				\addplot[line width=1pt,black,dotted]  table[x index=0, y index=2] {data/helm/sphere/sh0b_1_4_12};
			\end{semilogyaxis}
		\end{tikzpicture}
		\caption{$p=0$}
	\end{subfigure}
	\begin{subfigure}{.48\textwidth}
		\begin{tikzpicture}[scale = .95]
			\begin{semilogyaxis}[
					width=1\columnwidth,
					height=1\columnwidth,
					xlabel={$m$},
					ylabel={error},
					xticklabels={1,2,3,4},
					xtick = {1,2,3,4},
					grid=major,
				]
				\addplot[line width=1.5pt,mark = square,mark options={solid},mark size=3,red] table[x index=0, y index=3] {data/helm/sphere/sh1b_1_4_12};
				\addplot[line width=1.5pt,mark = square,mark options={solid},dotted,mark size=3,blue]  table[x index=0, y index=3] {data/helm/tours/sh1b_t_1_4_12};
				\addplot[line width=1.5pt,mark = square,mark options={solid},mark size=3]  table[x index=0, y index=3] {data/helm/ficher_data/sh1b_f_1_4_10};
				\addplot[line width=1,black,dotted] table[x index=0, y index=2] {data/helm/sphere/sh1b_1_4_12};
			\end{semilogyaxis}
		\end{tikzpicture}
		\caption{$p=1$}
	\end{subfigure}\\[.5cm]
	\begin{subfigure}{.48\textwidth}
		\begin{tikzpicture}[scale = .95]
			\begin{semilogyaxis}[
					width=1\columnwidth,
					height=1\columnwidth,
					xlabel={$m$},
					ylabel={error},
					xticklabels={1,2,3,4},
					xtick = {1,2,3,4},
					grid=major,
				]
				\addplot[line width=1.5pt,mark = square,mark options={solid},mark size=3,red]  table[x index=0, y index=3] {data/helm/sphere/sh2b_s_1_4_10};
				\addplot[line width=1.5pt,mark = square,mark options={solid},dotted,mark size=3,blue] table[x index=0, y index=3] {data/helm/tours/sh2b_t_1_4_12};
				\addplot[line width=1.5pt,mark = square,mark options={solid},mark size=3]  table[x index=0, y index=3] {data/helm/ficher_data/sh2b_f_1_4_10};
				\addplot[line width=1,black,dotted]  table[x index=0, y index=2] {data/helm/sphere/sh2b_s_1_4_10};
			\end{semilogyaxis}
		\end{tikzpicture}
		\caption{$p=2$}
	\end{subfigure}
	\begin{subfigure}{.48\textwidth}
		\begin{tikzpicture}[scale = .95]
			\begin{semilogyaxis}[
					width=1\columnwidth,
					height=1\columnwidth,
					xlabel={$m$},
					ylabel={error},
					xticklabels={1,2,3,4},
					xtick = {1,2,3,4},
					grid=major,
				]
				\addplot[line width=1.5pt,mark = square,mark options={solid},mark size=3,red]  table[x index=0, y index=3] {data/helm/sphere/sh3b_s_1_4_12};
				\addplot[line width=1.5pt,mark = square,mark options={solid},dotted,mark size=3,blue] table[x index=0, y index=3] {data/helm/tours/sh3b_t_1_4_12};
				\addplot[line width=1.5pt,mark = square,mark options={solid},mark size=3] table[x index=0, y index=3] {data/helm/ficher_data/sh3b_f_1_4_12};
				\addplot[line width=1,black,dotted]  table[x index=0, y index=2] {data/helm/sphere/sh3b_s_1_4_12};
			\end{semilogyaxis}
		\end{tikzpicture}
		\caption{$p=3$}
	\end{subfigure}

	\caption{Error for the Helmholtz problem with $\kappa = .5$. Maximal potential error. Sphere (red), Torus (blue, dotted) and Fichera cube (black).  The dotted line represents the expected potential convergence of $h^{2p+3}.$}
	\label{fig::numres::helm::sphere}
\end{figure}

In the case of the Helmholtz problem one observes, cf.~Figure \ref{fig::numres::helm::sphere}, that the convergence behaviour predicted by Corollary \ref{cor::pwconv} for smooth geometries (i.e.~sphere and torus) as well as for nonsmooth, nonconvex geometries (i.e.~the fichera cube), as shown up to order $p=3$.

For larger values of $p$, i.e.~$p>3$, the time required for the iterative solution of the system arising from \eqref{eq::varformdisc}, via a Jacobi preconditioned GMRES scheme, becomes prohibitive for refinement levels $m>3$. This can be attributed to the increased condition of the system due to the higher order approach. Thus, without the additional application of preconditioners tailored to the specific method, higher order boundary element methods for $p>3$ seem not to be advantageous over medium- and lower-order approaches, i.e.~$p\leq 3$.

\begin{figure}
	\centering
	\begin{subfigure}{.48\textwidth}
		\begin{tikzpicture}[scale = .95]
			\begin{semilogyaxis}[
					width=1\columnwidth,
					height=1\columnwidth,
					xlabel={$m$},
					ylabel={error},
					xticklabels={1,2,3,4},
					xtick = {1,2,3,4},
					grid=major,
				]
				\addplot[line width=1.5pt,mark = square,mark options={solid},mark size=3,red]  table[x index=0, y index=3] {data/helm/sphere/sh0b_1_4_12};
				\addplot[line width=1.5pt,mark = square,mark options={solid},dotted,mark size=3,blue]  table[x index=0, y index=3] {data/helm/sphere/sh0b_s2_1_4_10};
				\addplot[line width=1.5pt,mark = square,mark options={solid},mark size=3]  table[x index=0, y index=3] {data/helm/sphere/sh0b_s4_1_4_12};
				\addplot[line width=1,black,dotted]  table[x index=0, y index=2] {data/helm/sphere/sh0b_1_4_12};
			\end{semilogyaxis}
		\end{tikzpicture}
		\caption{$p=0$}
	\end{subfigure}
	\begin{subfigure}{.48\textwidth}
		\begin{tikzpicture}[scale = .95]
			\begin{semilogyaxis}[
					width=1\columnwidth,
					height=1\columnwidth,
					xlabel={$m$},
					ylabel={error},
					xticklabels={1,2,3,4},
					xtick = {1,2,3,4},
					grid=major,
				]
				\addplot[line width=1.5pt,mark = square,mark options={solid},mark size=3,red] table[x index=0, y index=3] {data/helm/sphere/sh1b_1_4_12};
				\addplot[line width=1.5pt,mark = square,mark options={solid},dotted,mark size=3,blue]  table[x index=0, y index=3] {data/helm/sphere/sh1b_s2_1_4_12};
				\addplot[line width=1.5pt,mark = square,mark options={solid},mark size=3]  table[x index=0, y index=3] {data/helm/sphere/sh1b_s4_1_4_12};
				\addplot[line width=1,black,dotted] table[x index=0, y index=2] {data/helm/sphere/sh1b_1_4_12};
			\end{semilogyaxis}
		\end{tikzpicture}
		\caption{$p=1$}
	\end{subfigure}\\[.5cm]
	\begin{subfigure}{.48\textwidth}
		\begin{tikzpicture}[scale = .95]
			\begin{semilogyaxis}[
					width=1\columnwidth,
					height=1\columnwidth,
					xlabel={$m$},
					ylabel={error},
					xticklabels={1,2,3,4},
					xtick = {1,2,3,4},
					grid=major,
				]
				\addplot[line width=1.5pt,mark = square,mark options={solid},mark size=3,red]  table[x index=0, y index=3] {data/helm/sphere/sh2b_s_1_4_10};
				\addplot[line width=1.5pt,mark = square,mark options={solid},dotted,mark size=3,blue] table[x index=0, y index=3] {data/helm/sphere/sh2b_s2_1_4_10};
				\addplot[line width=1.5pt,mark = square,mark options={solid},mark size=3]  table[x index=0, y index=3] {data/helm/sphere/sh2b_s4_1_4_12};
				\addplot[line width=1,black,dotted]  table[x index=0, y index=2] {data/helm/sphere/sh2b_s_1_4_10};
			\end{semilogyaxis}
		\end{tikzpicture}
		\caption{$p=2$}
	\end{subfigure}
	\begin{subfigure}{.48\textwidth}
		\begin{tikzpicture}[scale = .95]
			\begin{semilogyaxis}[
					width=1\columnwidth,
					height=1\columnwidth,
					xlabel={$m$},
					ylabel={error},
					xticklabels={1,2,3,4},
					xtick = {1,2,3,4},
					grid=major,
				]
				\addplot[line width=1.5pt,mark = square,mark options={solid},mark size=3,red]  table[x index=0, y index=3] {data/helm/sphere/sh3b_s_1_4_12};
				\addplot[line width=1.5pt,mark = square,mark options={solid},dotted,mark size=3,blue] table[x index=0, y index=3] {data/helm/sphere/sh3b_s2_1_4_10};
				\addplot[line width=1.5pt,mark = square,mark options={solid},mark size=3] table[x index=0, y index=3] {data/helm/sphere/sh3b_s4_1_4_10};
				\addplot[line width=1,black,dotted]  table[x index=0, y index=2] {data/helm/sphere/sh3b_s_1_4_12};
			\end{semilogyaxis}
		\end{tikzpicture}
		\caption{$p=3$}
	\end{subfigure}

	\caption{Error for the Helmholtz problem of the sphere. Maximal potential error. $\kappa=0.5$ (red), $\kappa = 2$ (blue, dotted) and $\kappa = 4$ (black).  The dotted line represents the expected potential convergence of $h^{2p+3}.$}
	\label{fig::numres::helm::kappa}
\end{figure}
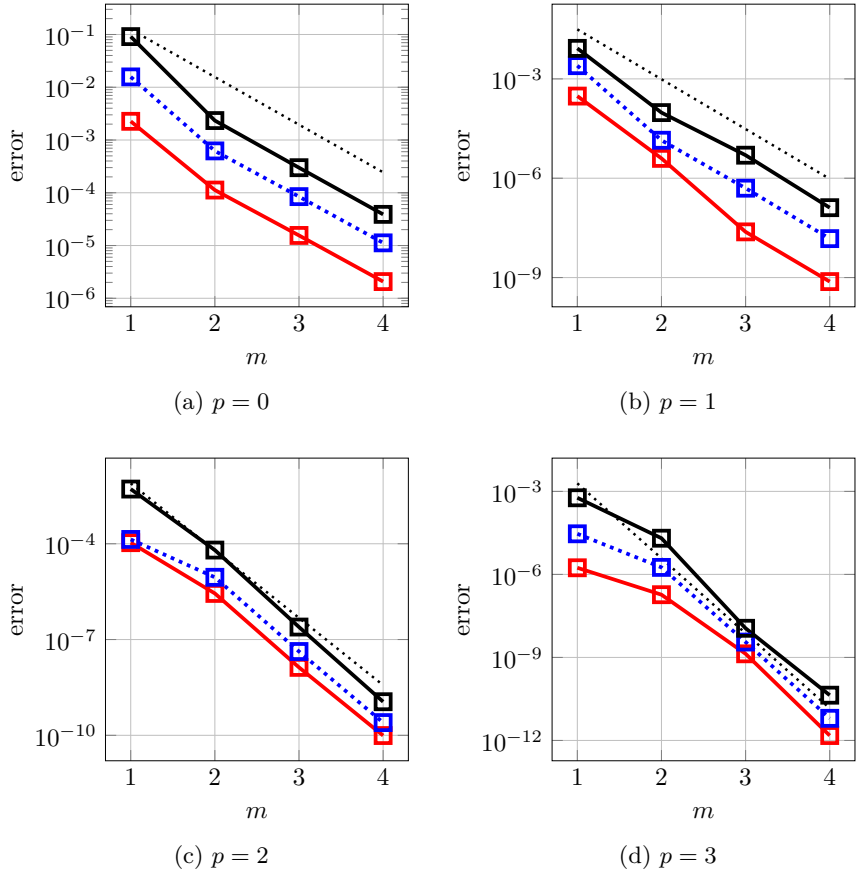

The convergence behaviour, as depicted in Figure \ref{fig::numres::helm::sphere} has also been observed for multiple, nonresonant wavenumbers, where the results for the case of the sphere and $\kappa = 0.5,1,2$ are showcased in Figure \ref{fig::numres::helm::kappa}. One can see nicely that $\kappa$ affects the constant appearing on the right hand side of estimate \eqref{eq::pwconv}, i.e.~for larger values of $\kappa$ the maximal potential error in the sense of \eqref{eq::pwconv} increases, as is to be expected of the Helmholtz problem. This can be seen since the wavenumber is encoded in the constant $C$ arising from $\norm{u^*(\bb x,\cdot)}_{H^{p+2}(\Gamma)}$ in \eqref{eq::constantcomesin}.

As in the Laplace case, orders of $p = 1,2,3$ prove advantageous in terms of time to solution compared to a lowest order approach with $p=0$, as can, once more, be seen in Table \ref{tab::time}.

\subsection{Time to Solution}\label{subsec::time}

The conjecture, that often higher order methods are, due to the high convergence rates, advantageous in terms of time to solution indeed goes well with the tendencies observed during the tests. However, since many of the computations have been done on different machines, partially with other computations in parallel, we chose to repeat some specific examples to make the times comparable with respect to order and refinement level.
Tables \ref{tab::time} and \ref{tab::time::fich} showcase said examples, where
$\operatorname{dof}^*$ gives the dimension of the corresponding space $\S^*_{p,m}(\Gamma)$, i.e.~corresponds to the size of the discrete system, and $\operatorname{dof}$ corresponds to the dimension of the corresponding spaces $\S_{p,m}(\Gamma)$, cf.\ Section \ref{subsec::embedding}.
\emph{Time} gives the time to solution in seconds, \emph{assembly} and \emph{solver} the corresponding times spent in matrix assembly and Jacobi preconditioned CG or GMRES, respectively.

All benchmarks were computed with the same codebase, described within Section \ref{sec::impl}.
The Laplace problem was computed on a Laptop, with an i7 7500U Kaby Lake processor (2 cores), 16GB RAM, and has been compiled with \texttt{g++ 5.4}, with compile flags \texttt{-O3 -march=native -fopenmp -flto}. The Helmholtz problem with $\kappa = 0.5$ was computed on a workstation, Xeon E5-2687W processor (4 cores), 16GB RAM, compiled via \texttt{g++ 4.8}, with compile flags \texttt{-O3 -march=native -fopenmp}.

However, as has been mentioned before, for $p>3$ and larger values of $\kappa,$ the time spent on solving the linear system increases, and eventually becomes prohibitive. Of course, the optimal choice of $p$ in times of solution depends on both, the accuracy requirements and the wavenumbers to consider.

\begin{table}
	\centering
	\caption{Exemplary time to solution for a lower and a higher order approach, for a maximal error of the potential of $\approx10^{-10}$ (Laplace) and $\approx10^{-9}$ (Helmholtz). $\operatorname{dof}^*=\operatorname{dim}\S^*_{p,m}(\Gamma)$, $\operatorname{dof}=\operatorname{dim}\S_{p,m}(\Gamma).$ For more detail see Section \ref{subsec::time}.}
	\label{tab::time}
	\begin{tabular}{c|l|r|r|c|r|r|r}\multicolumn{5}{c}{}\\
		    & \multicolumn{6}{c}{\textbf{Laplace problem, sphere geometry}}\\\hline
		    &            &        &        &     & \multicolumn{3}{c}{time (s)}\\
		$p$ & $m$ & accuracy   & $\operatorname{dof}^*$  & $\operatorname{dof}$      &  total    & assembly & solver   \\\hline
		$1$ & 6   & 4.0093e-11 & 98,304 & 25,350 &  2740.227 & 2683.600 & 56.627   \\\hline
		$4$ & 3   &  1.2673e-10 & 9,600  & 1,014  & 46.305   & 20.061   & 26.244   \\\hline \\
		    & \multicolumn{6}{c}{\textbf{Helmholtz problem, sphere geometry}}\\\hline
		    &            &        &        &     & \multicolumn{3}{c}{time (s)}\\
		$p$ & $m$ &accuracy   & $\operatorname{dof}^*$  & $\operatorname{dof}$    &  total    & assembly & solver   \\\hline
		$0$ & 7   &1.3521e-09 & 98,304 & 98,304 &  3106.417 & 3019.400 & 87.017  \\\hline
		$3$ &3   & 4.053e-09  & 6,144  & 726    &  273.741  & 2.951    & 270.790  \\\hline
	\end{tabular}
\end{table}



\begin{table}
	\centering
	\caption{Exemplary time to solution for a lower and a higher order approach, for a maximal error of the potential of $\approx10^{-5}$ (Laplace) and $\approx10^{-6}$ (Helmholtz). $\operatorname{dof}^*=\operatorname{dim}\S^*_{p,m}(\Gamma)$, $\operatorname{dof}=\operatorname{dim}\S_{p,m}(\Gamma).$ For more detail see Section \ref{subsec::time}.}
	\label{tab::time::fich}
	\begin{tabular}{c|l|r|r|c|r|r|r}\multicolumn{5}{c}{}\\
		    & \multicolumn{6}{c}{\textbf{Laplace problem, Fichera geometry}}\\\hline
		    &            &        &        &     & \multicolumn{3}{c}{time (s)}\\
		$p$ & $m$ & accuracy   & $\operatorname{dof}^*$  & $\operatorname{dof}$      &  total    & assembly & solver   \\\hline
		$0$ & 6 & 1.1397e-05 & 98,304 & 98,304 &  505.27 & 436.660 & 68.608 \\\hline
		$2$ & 3 & 1.2412e-05 & 13,824 & 2,400  &  49.690 & 10.102  & 39.588 \\\hline \\
		    & \multicolumn{6}{c}{\textbf{Helmholtz problem, Fichera geometry}}\\\hline
		    &            &        &        &     & \multicolumn{3}{c}{time (s)}\\
		$p$ & $m$ & accuracy   & $\operatorname{dof}^*$  & $\operatorname{dof}$    &  total    & assembly & solver   \\\hline
		$0$ & 4 & 1.3863e-06 & 24,576 & 24,576 &  56.562 & 45.821& 10.741\\\hline 
		$3$ & 2 & 6.9004e-07 & 6,144 & 1,176 &  13.473 & 0.803 & 12.67\\\hline
	\end{tabular}
\end{table}

\FloatBarrier
\section{Final Remarks}\label{sec::final}

The presented results show that higher order isogeometric boundary element methods achieve their corresponding orders of convergence, leading to an improvement in terms of time to solution. However, this effect cannot be expected for arbitrary orders and problems notorious for invoking ill conditioned systems. We will briefly summarise our findings below.

\begin{itemize}
	\item Indirect higher order boundary element methods play well with the core ideas of iso\-geo\-metric analysis. 
  The indirect approach is in\-sensitive w.r.t.~the choice of in- and exterior, and thus poses a black box scheme to solve PDEs in arbitrary areas. The high convergence rates promised by Corollary \ref{cor::pwconv} can indeed be reached in practice.
	\item In most test cases, higher order methods (with $p\leq 4$) are lucrative in terms of time to solution, even without nontrivial preconditioners. As seen in Tables \ref{tab::time} and \ref{tab::time::fich}, the application of higher order basis functions pays off, altough in case of the Fichera geometry one does not reach optimal convergence, cf.\ the regularity assumptions of Corollary \ref{cor::pwconv} and Figure \ref{fig::lapl::fichera}.
  While the iterative solvers take longer in the case of higher order ansatz functions due to ill conditioned systems, the time saved in matrix assembly proves to be significant for a sensible choice of $p$.
	\item The interpolation-based multipole method, in comparison with classical multipole methods, does not require a problem specific series expansion \cite{greengard2009fast,greengard1987fast}. It provides a general black box compression technique for boundary element methods, with milder requirements on the regularity of the kernel as other methods, cf.~Section \ref{subsec::fmm}.
	\item The approach via a discontinuous superspace, as explained in Section \ref{subsec::embedding}, offers enormous flexibility in the construction of different ansatz spaces. Moreover, the localization gained due to the artificial increase of dimension benefits an efficient implementation and the compression rates.
	\item The interpolation-based multipole approach might lead to instabilities in the kernel evaluation for $p\geq 5$ together with small $h$, since the order of the kernel interpolation in the near field requires high polynomial orders $\geq15$ for ansatz functions of order $4$ and up. This might lead to stability issues, if one were to take the order of ansatz functions higher than those presented in Section \ref{sec::numex}.
	\item Convergence rates of higher order approaches will suffer from insufficient regularity of geometry mappings (i.e.~due to knot repetition), even if the geometry might be smooth, cf.\ the remarks in the end of Section \ref{subsec::convergence}, Remark \ref{rem::smoothkernel} and Section \ref{subsec::perturbed}.
\end{itemize}

It should be noted that this is merely a first work into a direction that looks promising, since for medium orders of $p=2,3$ conditioning and issues of stability can still be handled by the usual means (i.e.,~Jacobi preconditioning).
Moreover, classical methods suffer from the the same problems as the method presented, regarding regularity and conditioning, as mentioned above \cite{Betcke_2011aa,Weggler_2011aa}.

\paragraph*{\textbf{Acknowledgement}} The work of J\"urgen D\"olz is supported by the Swiss National Science Foundation (SNSF) through the project \emph{$\mathcal H$-matrix based first and second moment analysis}. The work of Felix Wolf is supported by DFG Grants SCHO1562/3-1 and KU1553/4-1, the Excellence Initiative of the German Federal and State Governments and the Graduate School of Computational Engineering at TU Darmstadt.

\section*{References}
\bibliography{local}
\bibliographystyle{plain}
\end{document}